\documentclass [12pt]{article}
\usepackage{latexsym,amssymb,amsfonts,amsmath,amsthm,graphicx}
\pdfoutput=1
\usepackage{graphicx}
\usepackage{hyperref}

\newtheorem{lemma}{Lemma}[section]
\newtheorem{proposition}{Proposition}[section]
\newtheorem{corollary}{Corollary}[section]
\newtheorem{theorem}{Theorem}[section]
\newtheorem{remark}{Remark}[section]

\newtheorem{example}{Example}[section]

\newcommand{\s}{\vspace{2ex}}
\newcommand{\n}{\noindent}
\newcommand{\e}{\varepsilon}
\newcommand{\R}{\mathbb{R}}

\newcommand{\E}{\mathbb{E}}
\newcommand{\PP}{\mathbb{P}}

\title{Moment stability and large deviations for random dynamical systems on non-compact manifolds}

\author{Peter H. Baxendale} 

\AtEndDocument{\bigskip{\footnotesize
  \noindent
  \textsc{Peter H. Baxendale. Department of Mathematics, University of Southern California, CA 90089, USA.}
  \textit{E-mail address}: \texttt{\href{mailto:baxendal@usc.edu}{baxendal@usc.edu}}
 }}

\begin{document}

\maketitle

\begin{abstract}

Let $\{v_t: t \ge 0\}$, with values in $TM$, denote the linearization along trajectories of a stochastic differential equation on a Riemannian manifold $M$, and let $\lambda_t = t^{-1} \log \|v_t\|$ denote the associated finite time Lyapunov exponent.  
  The rate function for large deviations of $\lambda_t$ from its almost sure limit $\lambda$ is related via the G\"{a}rtner-Ellis theorem to the moment Lyapunov exponent $\Lambda(p) = \lim_{t \to \infty} t^{-1} \log \mathbb{E}[\|v_t\|^p]$.   
  When $M$ is compact there is a characterization of $\Lambda(p)$ in terms of an eigenvalue problem for an associated differential operator acting on functions on the unit sphere bundle $SM$.  We formulate growth conditions which, together with standard assumptions of hypoellipticity and positivity, ensure that the eigenvalue characterization of $\Lambda(p)$ remains valid on a suitable function space.  There is a central limit theorem for $\lambda_t$, and there are moderate deviation estimates for $\lambda_t$, both involving the second derivative $\Lambda''(0)$.  Examples are given showing the various different growth rates for $\Lambda(p)$ as $p \to \infty$ which may happen when $M$ is not compact.  

\end{abstract}

\s

\n Keywords: random dynamical system, Lyapunov exponent, moment stability, large deviations.

\section{Introduction} \label{sec intro} Consider the stochastic differential equation (SDE) for $\{x_t: t \ge 0\}$ on a $d$-dimensional $C^\infty$ manifold $M$ 
    \begin{equation} \label{x}
    dx_t = X_0(x_t)dt+ \sum_{j=1}^m X_j(x_t) \circ dW_t^j,
    \end{equation}
where the vector fields $X_0,X_1,\ldots,X_m$ are all $C^\infty$ and the $\{W_t^j: t \ge 0\}$, $j = 1, \ldots,m$, are independent standard Brownian motions, and write
   \begin{equation} \label{L}
     L = \frac{1}{2}\sum_{j=1}^m X_j^2 + X_0.
     \end{equation}
Consider also the real valued additive functional
   \begin{align}      
      A_t & = \int_0^t q_0(x_s)ds + \sum_{j=1}^m \int_0^t q_j(x_s) \circ dW_s^j   =  \int_0^t Q(x_s)ds + \sum_{j=1}^m \int_0^t q_j(x_s) dW_s^j \label{A}
    \end{align}     
where $q_0, q_1, \ldots,q_m$ are all in $C^\infty(M)$ and $Q(x) = q_0(x)+ \frac{1}{2} \sum_{j=1}^m X_jq_j(x)$.    
Our interest is in the behavior of $A_t$ as $t \to\infty$.  In particular we study the logarithmic moment generating function
    \begin{equation} \label{Lam}
      \Lambda(p) = \lim_{t \to\infty} \frac{1}{t} \log \E^x\left[ e^{pA_t}\right], \quad p \in \R,
    \end{equation}      
and especially its characterization in terms of an eigenvalue problem on a suitable function space.

This problem arises when studying the stability properties of the solution $\{v_t: t \ge 0\}$ of the linear SDE in $\R^d$
  \begin{equation} \label{v}
    dv_t = B_0 v_t dt+ \sum_{j=1}^mB_j v_t \circ dW_t^j, \quad v_0 \neq 0.
    \end{equation} 
Following Khas'minskii \cite{Khas67}, see also \cite[Section 6.7]{KhAS}, write $\theta_t = v_t/\|v_t\|$.  Then $\{\theta_t: t \ge 0\}$ is a diffusion process in $\mathbb{S}^{d-1}$ given by
  \begin{equation} \label{theta}
    d\theta_t = \bigl(B_0 \theta_t - \langle B_0 \theta_t,\theta_t \rangle \theta_t\bigr) dt+ \sum_{j=1}^m \bigl(B_j \theta_t - \langle B_j \theta_t,\theta_t \rangle \theta_t\bigr) \circ dW_t^j,
    \end{equation} 
and 
   \begin{equation} \label{logv}
   \log \|v_t\| - \log \| v_0\| = \int_0^t \langle B_0 \theta_s,\theta_s\rangle ds + \sum_{j=1}^m \int_0^t \langle B_j \theta_s,\theta_s\rangle \circ dW_s^j \equiv A_t.
  \end{equation}
With $\{\theta_t: t \ge 0\}$ on $\mathbb{S}^{d-1}$ given by \eqref{theta} and the addititive functional $\{A_t: t \ge 0\}$ given by \eqref{logv}, the function $\Lambda(p)$ can now be recognized as the $p$th moment Lyapunov exponent
   $$
   \Lambda(p) = \lim_{t \to \infty}\frac{1}{t} \log \E^\theta\left[e^{p \log \|v_t\|}\right] = \lim_{t \to \infty}\frac{1}{t} \log \E^\theta\left[ \|v_t\|^p \right].
   $$
A similar situation arises with the random linear ODE in $\R^n$
    \begin{equation} \label{v2}
      dv_t = B(x_t)v_tdt
    \end{equation}  
where $\{x_t: t \ge 0\}$ is the diffusion on $M$ given by \eqref{x} and $B$ is a mapping from $M$ into the space of $n \times n$ real matrices.  In this case, writing $\theta_t = v_t/\|v_t\| \in \mathbb{S}^{n-1}$ gives
     \begin{equation} \label{theta2}
    d\theta_t = \bigl(B(x_t) \theta_t - \langle B(x_t) \theta_t,\theta_t \rangle \theta_t\bigr) dt,
    \end{equation} 
and 
   \begin{equation} \label{logv2}
   \log \|v_t\| - \log \| v_0\| = \int_0^t \langle B(x_s) \theta_s,\theta_s\rangle ds \equiv A_t. 
  \end{equation} 
The pair of equations (\ref{x},\ref{theta2}) gives an SDE for the diffusion $\{(x_t,\theta_t): t \ge 0\}$ on $M \times \mathbb{S}^{n-1}$, and then \eqref{logv2} gives the functional $\{A_t: t \ge 0\}$ with $Q(x,\theta) = \langle B(x)\theta,\theta \rangle: M \times \mathbb{S}^{n-1} \to \R$.

The relationship between almost sure and moment stability for the random linear system (\ref{x},\ref{v2}) in the case of compact $M$ is discussed using the function $\Lambda(p)$ by Molchanov \cite{Mol78} and Arnold \cite{Arn84}.  See also Arnold, Kliemann and Oeljeklaus \cite{AKO86}.  Results for the linear SDE \eqref{v} are given in Arnold, Oeljeklaus and Pardoux \cite{AOP86}.  For connections with the large deviation behavior of $\frac{1}{t} \log \|v_t\|$ as $t \to \infty$ and an associated central limit theorem see Stroock \cite{Str86}, Baxendale \cite{Bax87}, and Arnold and Kliemann \cite{AK87}.

More generally, given the (possibly nonlinear) SDE \eqref{x} it is natural to consider the dependence of the solution $\{x_t: t \ge 0\}$ upon its initial value $x_0$.  The effect of an infinitesimal change in the initial condition can be seen in the behavior of the linearized system $\{v_t: t \ge 0\}$ with values in the tangent bundle $TM$ given by 
   \begin{equation} \label{Tx}
    dv_t = TX_0(v_t)dt+ \sum_{j=1}^m TX_j(v_t) \circ dW_t^j, \quad v_0 \neq 0
    \end{equation}
where $TX_j$ denotes the derivative of $X_j$, thought of as a map $TM \to T(TM)$.  Assuming the manifold $M$ has a Riemannian structure, let $\|v\|$ denote the Riemannian norm of the tangent vector $v \in TM$ and let $SM = \{v \in TM: \|v\| = 1\}$ denote the unit sphere bundle.    
   Because $TX_j$ is linear on each tangent space $T_xM$ the process $\{\theta_t = v_t/\|v_t\|: t \ge 0\}$ is a diffusion in $SM$.  From \eqref{Tx} we obtain 
    \begin{equation} \label{Sx}
    d\theta_t = \widetilde{X}_0(\theta_t)dt+ \sum_{j=1}^m \widetilde{X}_j(\theta_t) \circ dW_t^j, \quad \theta_0 = v_0/\|v_0\|
    \end{equation}
where $\widetilde{X}_j(\theta)$ is the orthogonal projection of $TX_j(\theta)$ onto $T_\theta SM$.  In terms of the decomposition of $T_vTM$ 
determined by the Levi-Civita connection $\nabla$ associated with the Riemannian structure on $M$, the horizontal and vertical components of $\widetilde{X}_j(\theta)$ for $\theta \in S_xM$ are $X_j(x)$ and $\nabla X_j(x) \theta - \langle \nabla X_j(x) \theta,\theta \rangle \theta$ respectively. 
Also from \eqref{Tx} we have 
   \begin{equation} \label{logTx}
   \log \|v_t\| - \log \| v_0\| = \int_0^t q_0(\theta_s) ds + \sum_{j=1}^m \int_0^t q_j(\theta_s) \circ dW_s^j \equiv A_t,
  \end{equation}
where $q_j(\theta) = \langle \nabla X_j(x) \theta,\theta \rangle $ for $\theta \in S_xM$.   See \cite[page 189]{BS88} for details of this calculation.  
The pairing of the diffusion $\{\theta_t: t \ge 0\}$ on $SM$ given by \eqref{Sx} together with the additive functional $\{A_t: t \ge 0\}$ given by \eqref{logTx} again fits into the general setting of (\ref{x},\ref{A}).

Results for the case when $M$ is compact were obtained in Baxendale and Stroock \cite{BS88}.  They include the existence and characterization of $\Lambda(p)$ in terms of an eigenvalue 
problem, and also large deviation estimates and a central limit theorem for the limiting behavior of $A_t/t$ as $t \to \infty$.  In the case where the top Lyapunov exponent is positive, the eigenfunctions associated with $\Lambda(p)$ were used in \cite{BS88} to obtain rigorous estimates on exponential moments of exit times for the 2-point motion 
  $\{(x_t,y_t): t \ge 0\}  $ 
from a neighborhood of the diagonal in $M\times M$, leading to almost sure mixing results for the associated stochastic flow of diffeomorphisms of $M$, see Dolgopyat, Kaloshin \& Koralov \cite{DKK}.  

\s

The aim of the present paper is to obtain corresponding results for the case of a noncompact manifold $M$.   It is enough to study the original system (\ref{x},\ref{A}) on $M$, as any results proved in that setting will apply also when $M$ is replaced by $M \times \mathbb{S}^{n-1}$ or $SM$.  In order to obtain results about $\Lambda(p)$, we first study the richer structure of the twisted semigroup (or Feynman-Kac semigroup) of operators $\{T_t^p: t \ge 0\}$ acting on some space of functions $g: M \to \R$ by the formula
   \begin{equation} \label{Ttp}
   T_t^p g(x) = \E^x\left[g(x_t)e^{p A_t}\right], \quad x \in M.
  \end{equation}
Information about the operators $T_t^p$ is then used to evaluate and characterize
     $$
     \Lambda(p) = \lim_{t \to \infty} \frac{1}{t} \log T_t^p \mathbf{1}(x).
     $$  
In order for this program to succeed it is necessary to find a suitable space of functions for the $T_t^p$ to operate on.  When $M$ is compact it is natural to consider the space of bounded measurable functions, or the space of bounded continuous functions.  See for example \cite{Bax87} for the linear SDE \eqref{v}.  But for the case where $M$ is non-compact we will need to impose certain growth conditions on the operator $L$ and the functions $Q$ and $q_1, \ldots,q_m$. 

There is much related work, mostly in the case when all the $q_j$ are zero, because of the relationship of $\Lambda(p)$ via the G\"{a}rtner-Ellis theorem to the large deviation theory for $\frac{1}{t} \int_0^t Q(x_s)ds$  and more generally to the large deviation theory for the empirical measure $\frac{1}{t}\int_0^t \delta_{x_s}ds$.  Note that the large deviation behavior for the empirical measure requires the consideration of a large class of functions $Q$, including non-smooth functions, while we will restrict to a single given smooth function $Q$.

 Especially relevant for our purposes are the papers by Liming Wu \cite{WuL95}, \cite{WuL01}, Kontoyiannis and Meyn \cite{KM03}, \cite{KM05}, Ferr\'{e} and Stoltz \cite{FS20}, and Ferr\'{e}, Rousset and Stoltz \cite{FRS21}.   
 Following these papers, given a function $V: M \to [1,\infty)$, we introduce the real Banach space 
  \begin{equation} \label{BV}
  B_V = \Bigl\{\mbox{measurable }g: M \to \R: \|g\|_V \equiv \sup\{|g(x)|/V(x): x \in M\} < \infty\Bigr\}.
  \end{equation}
Assuming $V$ satisfies certain Lyapunov type conditions, see the growth assumption {\bf (G)} in Section \ref{sec ass}, then $\{T_t^p: t \ge 0\}$ acts as a semigroup of bounded operators on $B_V$.  Under additional conditions that there is enough noise in the SDE \eqref{x}, see the hypoellipticity and positivity assumptions {\bf(H)} and {\bf (P)} in Section \ref{sec ass}, we show that $\Lambda(p)$ is well-defined and is characterized by the eigenvalue/eigenfunction equation 
   \begin{equation} \label{eigen}
   \Bigl(L+ p \sum_{j=1}^m q_j X_j + pQ + \frac{1}{2}p^2 \sum_{j=1}^m q_j^2\Bigr) g(x) = \Lambda(p)g(x)
  \end{equation}
for positive $g \in B_V$, see Corollary \ref{cor Lam}. We note that this characterization may fail when the requirement that $g \in B_V$ is removed, see Example \ref{ex OU}. 

The assumptions {\bf (H)}, {\bf (P)} and {\bf(G)} are given in Section \ref{sec ass} and the main results are given in Section \ref{sec results}.  
In addition to the characterization of $\Lambda(p)$ there are results on the real analyticity of the mappings $p \to \Lambda(p)$ and $p \to \phi_p \in B_V$, where $\phi_p$ denotes the eigenfunction in \eqref{eigen}.   
Properties of $\Lambda(p)$ and the corresponding eigenfunctions $\phi_p$ for $p$ near zero are used to describe the almost sure limit of $A_t/t$ and the variance in a central limit theorem for $A_t/t$ in terms of $\Lambda'(0)$ and $\Lambda''(0)$.  
Section \ref{sec examples} contains two examples where these theorems can apply to issues of moment stability for tangent vectors, and Section \ref{sec asymp} contains examples showing how the growth of $\Lambda(p)$ as $p \to \infty$ when $M$ is not compact can be very different from the behavior in the compact case.  
Finally, Sections \ref{sec fix proof}, \ref{sec vary proof} and \ref{sec  near proof} contain the proofs for the main results in Section \ref{sec results}, and Section \ref{sec asymp proof} contains the calculations for the examples in Section \ref{sec asymp}.

\section{Assumptions} \label{sec ass}

Throughout $M$ is a $C^\infty$ manifold of dimension $d$.  The diffusion $\{x_t: t \ge 0\}$ on $M$ is given by the SDE \eqref{x}
and the additive functional $\{A_t: t \ge 0\}$ is given by \eqref{A}.
The vector fields $X_0, \ldots, X_m$ and the functions $Q, q_0, q_1, \ldots,q_m$ are all assumed to be $C^\infty$. 
The second order linear differential operator $L$ is given by \eqref{L} and its action 
on $C^2$ functions with compact support agrees with the action of the extended generator (see \cite{MTIII}) of the Markov process $\{x_t: t \ge 0\}$.      

We will make three assumptions about the system (\ref{x},\ref{A}).  The first two assumptions {\bf (H)} and {\bf (P)} will ensure that the noise in the SDE \eqref{x} is sufficient to cause smoothing and positivity respectively.  They will be satisfied if $L$ is elliptic.  The third assumption involves a Lyapunov type function $V: M \to [1,\infty)$ and will be used to control the behavior of $x_t$ and $A_t$ outside compact subsets of $M$. 

Assumptions similar to {\bf(H)} and {\bf (P)} appear in \cite{FS20}.  Other papers replace {\bf (H)} and {\bf (P)} with assumptions of smoothness and positivity of the integral kernel of $T_t^p$, see \cite[Assumption 5]{FRS21}.  
In the case with $q_j \equiv 0$ the growth condition {\bf (G)} is very similar to \cite[Assumption 4]{FRS21}.

\subsection{Hypoellipticity}  Let ${\cal L} = {\cal L}(X_0,X_1,\ldots,X_m)$ be the Lie algebra generated by the vector fields $X_0,X_1,\ldots,X_m$ and then let ${\cal L}_0 = {\cal L}_0(X_0;X_1,\ldots,X_m)$ be the ideal in ${\cal L}$ generated by $X_1,\ldots,X_m$.  

\s

\n{\bf Assumption (H).}   Assume  $\mbox{span}\{U(x): U \in {\cal L}_0\} = T_xM$ for all $x \in M$.     
  
\s

The assumption {\bf (H)} is often referred to as the parabolic  H\"{o}rmander condition, see for example Hairer \cite{Hai}.  
It 
is equivalent to either of the statements that 
       $$
      \mbox{span}\left\{U(t,x): U \in {\cal L}\left(X_0 - \frac{\partial}{\partial t}, X_1,\ldots,X_m\right)\right\} = T_{(t,x)}(\R \times M) \quad \mbox{ for all } (t,x) \in \R \times M
       $$
or 
   $$
      \mbox{span}\left\{U(t,x): U \in {\cal L}\left(X_0 + \frac{\partial}{\partial t}, X_1,\ldots,X_m\right)\right\} = T_{(t,x)}(\R \times M) \quad \mbox{ for all } (t,x) \in \R \times M.
       $$         
Then by H\"{o}rmander \cite{hor} the assumption {\bf (H)} implies $L  - \frac{\partial}{\partial t}$ and $L  + \frac{\partial}{\partial t}$ are hypoelliptic.


\subsection{Positivity}  

Let $\PP_t(x,A) = \PP^x(x_t \in A)$ denote the transition probability for the diffusion $\{x_t: t \ge 0\}$. 

\s

\n{\bf Assumption (P).}  $\PP_t(x,U) > 0$ for all $t > 0$, $x \in M$ and non-empty $U \subset M$.

\s

This assumption can be verified using the following controllability condition.   For $t>0$ let ${\cal U}_{t}$ denote the space of piecewise continuous paths $u:[0,t] \to \R^m$, and for $u \in {\cal U}_{t}$ let   $\{\xi(s,x;u):0 \leq s\leq t\}$ denote the controlled path 
$$
\frac{\partial \xi}{\partial s}(s,x;u) = X_0(\xi(s,x;u))  +   \sum_{j=1}^m u_j(s)X_j(\xi(s,x;u))
$$
with $\xi(0,x;u)=x$.

\s

\n{\bf Assumption (C).}  The set $\{\xi(t,x;u): u \in {\cal U}_t\}$ is dense in $M$ for all $t > 0$ and $x \in M$.  

\s

The Stroock-Varadhan support theorem \cite{SVsupp} shows that {\bf (C)} implies {\bf (P)}.


\subsection{Growth} 

  Define the vector field $Y(x)= \sum_{j=1}^m q_j(x)X_j(x)$ and the function $R(x) = \sum_{j=1}^m \bigl(q_j(x)\bigr)^2$.  For $p \in\R$ define the second order linear differential operator $L_p$ acting on functions $g:M \to \R$ by 
    \begin{equation} \label{Lp}
    L_pg(x) = Lg(x) + pYg(x) + pQ(x)g(x) + \frac{p^2}{2}R(x)g(x). 
        \end{equation}
That is, $L_p = L+pY+pQ+ \frac{p^2}{2}R$.  The assumption below is for a fixed value of $p$.

\s

\n{\bf Assumption (G).}  There exists a $C^2$ function $V:M  \to [1,\infty)$ with $V(x) \to \infty$ as $x \to \infty$ so that 
\begin{description}
     \item[(0)] $\dfrac{LV(x)}{V(x)}$ is bounded above. 
    \item[(1)] $\dfrac{L_pV(x)}{V(x)}  \to -\infty$ as $x \to \infty$. 
     \item[(2)] There exists $\beta > 1$ such that $\dfrac{L_{\beta p}V(x)}{V(x)}$ is bounded above.    
    \item[(3)] There exists $\beta > 1$ such that $\dfrac{L_p(V^\beta)(x)}{V^\beta(x)}$ is bounded above.    
     \end{description}
Note $\beta$ can be different in {\bf (2)} and {\bf (3)}.   

\begin{remark} \label{rem noncompact} This assumption is designed for the case when $M$ is noncompact, and then ``$V(x) \to \infty$ as $x \to \infty$'' is shorthand for ``$\{x \in M: V(x) \le k\}$ is compact for all $k \in \R$''.  Similarly for {\bf (1)}.  Of course, if $M$ happens to be compact then {\bf (G)} is trivially true with $V(x) \equiv 1$.  
 \end{remark}

\begin{remark} Write $V(x) = e^{\gamma U(x)}$ with $U(x) \ge 0$ and $U(x) \to \infty$ as $x \to \infty$.  The left side of {\bf (1)} is 
\begin{equation} \label{LpV}
\frac{L_pV(x)}{V(x)} = \gamma LU(x)+ \frac{\gamma^2}{2}\sum_{j=1}^m \bigl((X_j U)(x)\bigr)^2 +p \gamma YU(x)+pQ(x)+ \frac{p^2}{2}R(x).
\end{equation}
Similarly for {\bf (2)} and {\bf (3)}.  The convexity of the right side of \eqref{LpV} with respect to $p$ and with respect to $\gamma$ implies that if {\bf (1)} is satisfied then $\beta$ in {\bf (2)} or {\bf (3)} can be replaced by any smaller $\tilde{\beta} >1$.  It follows that without loss of generality we may assume the same value of $\beta$ in {\bf (2)} and {\bf (3)}. \end{remark}

\begin{remark} \label{rem explosion}  The condition {\bf(0)} is included to ensure that the solution $\{x_t: t \ge 0\}$ of \eqref{x} exists for all $t \ge 0$.  The condition {\bf (1)} is designed to ensure that the $\{T_t^p: t \ge 0\}$ act as bounded operators on the Banch space $B_V$.  The conditions {\bf (2)} and {\bf (3)} are more technical and are used to prove uniform integrability at certain points in the proofs.

\end{remark}

\section{Statement of results} \label{sec results}

This section is divided into three subsections.  Section \ref{sec fix} deals with $\Lambda(p)$ for a single fixed $p$ and includes the important characterization of $\Lambda(p)$ as the solution of an eigenvalue/eigenfunction problem on the Banach space $B_V$, see Corollary \ref{cor Lam}.  Section \ref{sec vary} deals with the dependence of $\Lambda(p)$ and the associated eigenfunction $\phi_p$ on the parameter $p$.  Section \ref{sec near} deals with the behavior of $\Lambda(p)$ for $p$ near 0.  It includes the identification of $\Lambda'(0)$ as the almost sure limit of $A_t/t$ as $t \to \infty$, and of $\Lambda''(0)$ as the variance in a central limit theorem for $A_t/t$ as $t \to \infty$, see Theorem \ref{thm clt}.

   \subsection{Results for fixed $p$} \label{sec fix}

Assume {\bf (G)} with the function $V$, and let $B_V$ be the Banach space defined in \eqref{BV}.  For linear $T: B_V \to B_V$ write the operator norm $\|T\|_V = \sup\{ \|Tg\|_V: \|g\|_V = 1\}$.  Recall from \eqref{Ttp} the definition
  \begin{equation} \nonumber
   T_t^p g(x) = \E^x\left[g(x_t)e^{p A_t}\right]
  \end{equation}
for measurable $g: M \to \R$.  The first result does not require the positivity assumption {\bf (P)}.   The main novelty here is part (ii) which relies on the smoothness of $Q$ and the $q_j$.   In the case when $q_j \equiv 0$, parts (i) and (iii) are valid under similar growth conditions for a wider class of functions $Q$,  see \cite[Lemma 2.6]{FS20} and \cite[Theorem 3]{FRS21}.

\begin{theorem}  \label{thm Tt} Assume {\bf (H,G)}.

(i)  $\{T_t^p: t \ge 0\}$ acts as a semigroup of bounded operators on $B_V$.

(ii) For $g \in B_V$ the mapping $(t,x) \to T_t^pg(x)$ is $C^\infty$ on $(0,\infty) \times M$ and 
     $$
     \left(L+pY+pQ+\frac{p^2}{2}R - \frac{\partial}{\partial t}\right) T_t^pg(x) = 0.
     $$

(iii)  For $t > 0$ the mapping $T_t^p$ is a compact operator on $B_V$.

\end{theorem}

The next result applies the Kre\u{\i}n-Rutman theorem for positive compact linear operators \cite{KR50} to the operator $T_1^p$ to obtain the eigenfunction $\phi_p$, and then applies the theory of $\widehat{V}$-exponential ergodicity \cite{MTbook} to the twisted Markov operator $\widehat{T}g(x) =T_1^p(\phi_p g)(x)/ \bigl(\kappa_p \phi_p(x)\bigr)$ where $\widehat{V}(x) = V(x)/\phi_p(x)$.  This is the method used in \cite{KM05} and \cite{FRS21}.  For completeness we give a precise statement here and the proof later. 

\begin{theorem} \label{thm KR} Assume {\bf (H,P,G)}.   Let $\kappa_p$ be the spectral radius of the operator $T_1^p$.   Then $\kappa_p$ is a positive eigenvalue of $T_1^p$ with corresponding eigenfunction $\phi_p \in B_V$ with $\phi_p(x) > 0$ for all $x$.  There exists a measure $\widehat{\nu}_p$ with smooth density on $M$ such that $\int_M V d\widehat{\nu}_p < \infty$ and $\widehat{\nu}_p(U) > 0$ for every non-empty open $U \subset M$, and there exist constants $C_p < \infty$ and $0 < \eta_p< 1$ such that if $g \in B_V$ then
  \begin{equation} \label{get}
  \left|\frac{1}{(\kappa_p)^t} T_t^p g(x) - \phi_p(x) \int_M g d\widehat{\nu}_p \right|  \le C_p \|g\|_V \eta_p^t V(x)
  \end{equation}
for all $x \in M$ and $t >0$.  Moreover $\phi_p$ is an eigenfunction of $T_t^p$ with eigenvalue $(\kappa_p)^t$.  Finally, if $g \in \mathbb{C}B_V$ is a (non-trivial) eigenfunction of the complexification of $T_t^p$ for some $t >0$ then either $\int g\,d\widehat{\nu}_p = 0$ and the eigenvalue has modulus strictly less than $(\kappa_p)^t$, or else $g$ is a (complex) multiple of $\phi_p$.
  \end{theorem}

With these results about the operators $T_t^p$ we can take $g$ to be identically 1 and obtain results about $\E^x\left[e^{pA_t}\right] = T_t^p\mathbf{1}(x)$.

\begin{remark}  In the proof of Theorem \ref{thm KR}, once we have the eigenfunction $\phi_p$ then the measure $\widehat{\nu}_p$  is determined by the condition $\int \phi_p d\widehat{\nu}_p = 1$.  Conversely, we can scale the finite measure $\widehat{\nu}_p$ to be a probability measure and then scale the eigenfunction to satisfy $\int \phi_p d\widehat{\nu}_p = 1$.  In this case, taking $g = \mathbf{1}$ in \eqref{get} we get
      $$
      \left| \frac{1}{ (\kappa_p)^t}\E^x\left[e^{p A_t}\right] - \phi_p(x) \right| \le C_p\eta_p^t V(x).
      $$
Kontoyiannis and Meyn refer to this as a multiplicative ergodic theorem, see \cite[Theorem 3.1(iii)]{KM05}.

\end{remark}

\begin{corollary} \label{cor Lam}  Assume {\bf (H,P,G)}.

(i)  For all $x \in M$
   $$
   \lim_{t \to \infty} \frac{1}{t} \log \E^x\left[e^{p A_t}\right] = \log \kappa_p,
   $$
so that $\Lambda(p)$ is well-defined in \eqref{Lam} and $\Lambda(p) = \log \kappa_p$.

(ii)  The positive eigenfunction $\phi_p$ is in $B_V \cap C^\infty(M)$ and satisfies 
      $$
      \left(L+pY+p Q+\frac{p^2}{2}R\right)\phi_p = \Lambda(p) \phi_p.
      $$
      
(iii) Suppose $g \in B_V \cap C^2(M)$ is non-negative, not identically 0, and satisfies 
     $$
      \left(L+pY+p Q+\frac{p^2}{2}R\right)g = \widetilde{\Lambda} g.
      $$
      for some $\widetilde{\Lambda}$.  Then $\widetilde{\Lambda} = \Lambda(p)$ and $g$ is a positive multiple of $\phi_p$.    
       
\end{corollary}       

The proofs of Theorems \ref{thm Tt} and \ref{thm KR} and Corollary \ref{cor Lam} are given in Section \ref{sec fix proof}. 

\s

Often there is a choice to be made in the function $V$.  Clearly the existence of $\Lambda(p)$ in (i) is not affected by the choice of $V$.  But for (ii) smaller $V$ is more useful because more can be said about the eigenfunction $\phi_p$, and for (iii) larger $V$ is more useful because it makes it easier to check that a candidate function $g$ is an element of $B_V$.  

\begin{example} \label{ex pitch} {\rm The pitchfork bifurcation with additive noise
     $$
     dx_t = (ax_t - b x_t^3)dt+\sigma dW_t,
     $$
see \cite{CF98}, \cite{CDLR}, \cite{BBBE}, has linearization
     $$
     dv_t = (a-3bx_t^2)v_t dt.
     $$
The $p$th moment Lyapunov exponent $\Lambda(p)= \lim_{t \to \infty} \frac{1}{t} \log \E^x[|v_t|^p]$ corresponds to $M = \R$ and $L = \frac{1}{2} \sigma^2 \frac{d^2}{dx^2} + (ax-bx^3)\frac{d}{dx}$ and $Q(x) = a-3x^2$ and $q(x) = 0$.  Since $L$ is elliptic the conditions {\bf (H)} and {\bf (P)} are satisfied.  It is easy to verify {\bf (G)} when $V$ is taken to be any one of the following:
   \begin{itemize}
    \item $V(x) = (1+x^2)^\alpha$ for any $\alpha > \max(0,-3p/2)$,
    \item $V(x) = e^{\gamma x^2}$ for any $\gamma > 0$,
    \item $V(x) = e^{\gamma x^4}$ for $0 < \gamma < b/(2 \sigma^2)$.
      \end{itemize}    
}\end{example}

\s

The following cautionary example shows the importance of part (iii) of Corollary \ref{cor Lam}.

\begin{example} \label{ex OU} {\rm Consider the Ornstein-Uhlenbeck process 
    $$
    dx_t = -a x_t dt + \sigma dW_t
    $$
(with $a > 0$) and $A_t = \int_0^t x_s^2 ds$.  This is an example of a quadratic functional of a Gaussian process, see Bryc and Dembo \cite{BD97}.   Fix $p < a^2/(2\sigma^2)$.  The eigenvalue equation 
    $$
    (L+px^2) g(x) = \widetilde{\Lambda}g(x)
    $$
has infinitely many distinct solutions with positive $g$, including the solutions
    $$
    (L+px^2)e^{\beta_1 x^2} =  \sigma^2 \beta_1 e^{\beta_1 x^2}
    $$
and 
   $$
    (L+px^2)e^{\beta_2 x^2} =  \sigma^2 \beta_2 e^{\beta_2 x^2},
    $$  
where $\beta_1$ and $\beta_2$ are given by 
   $$
    \beta_1 = \frac{a - \sqrt{a^2-2 p \sigma^2}}{2\sigma^2} < \frac{a}{2\sigma^2} < \frac{a + \sqrt{a^2-2 p \sigma^2}}{2\sigma^2} = \beta_2.
  $$
If we knew only of the second solution, we might be tempted to claim $\Lambda(p) = \sigma^2\beta_2$.  
However, taking $V(x) = e^{\gamma x^2}$ with $\gamma = a/(2\sigma^2)$, we have 
  $$
   \frac{(L+px^2)V(x)}{V(x)} = \gamma \sigma^2 + \Bigl( -2 \gamma a + 2 \gamma^2 \sigma^2 + p \Bigr)x^2= \frac{a}{2} + \left(p-\frac{a^2}{2\sigma^2}\right)x^2.
   $$
 For $p < a^2/(2\sigma^2)$ the condition {\bf (G)} is satisfied with $V(x) = e^{ax^2/(2 \sigma^2)}$, and we may apply Corollary \ref{cor Lam}(iii). 
  Since $\beta_1 < a/(2\sigma^2)$ the function $e^{\beta_1 x^2}$ is in $B_V$ and so  
    $$
    \Lambda(p) = \sigma^2 \beta_1 = \frac{1}{2}(a-\sqrt{a^2-2p \sigma^2}), \qquad p < \frac{a^2}{2\sigma^2}.
    $$
Note that for $p \ge a^2/(2 \sigma^2)$ the condition {\bf (G)} is not satisfied with $V(x)$ of the form $e^{\gamma x^2}$ for any $\gamma$.  Alternative methods, see \cite{BD97}, give $\Lambda(a^2/(2\sigma^2)) = a/(2\sigma^2)$ and $\Lambda(p) = \infty$ for $p > a^2/(2 \sigma^2)$.

}\end{example}

\subsection{Results for varying $p$} \label{sec vary}

Suppose that the growth condition {\bf (G)} is satisfied at two values $p_1 < p_2$ with the same function $V$.  The convexity of \eqref{LpV} with respect to $p$ implies that {\bf (G)} is satisfied for all $p \in [p_1,p_2]$.    Assuming {\bf (H)} and {\bf (P)} are true for the diffusion $\{x_t: t \ge 0\}$, then $\Lambda(p)$ is well-defined and finite for all $p \in [p_1,p_2]$.  

The results in this Section use the analytic perturbation theory of Kato \cite{Kato} for operators acting on a complex vector space.  The crucial fact is that the operators $\{T_t^p: p \in (p_1,p_2)\}$ acting on $B_V$ for real $p$ can be extended to an analytic family of operators $T_t^z$ acting on the complexification $\mathbb{C}B_V$ of $B_V$ as $z$ runs through some neighborhood of $(p_1,p_2)$ inside $\mathbb{C}$, see Lemma \ref{lem T anal}.   This observation is well known for the case $q_j \equiv 0$ under various growth conditions, see for example \cite[Theorem 3.1]{KM05}, \cite[Lemma 6.7]{FS20} and \cite[Theorem 2.10]{BBBE}.

\begin{theorem} \label{thm kato1}   Assume {\bf (H,P,G)} valid for two values $p_1 < p_2$ with same function $V$.  Then $p \to \Lambda(p)$ is convex and real analytic on $(p_1,p_2)$ and there is a real analytic mapping $p \to \phi_p$ from $(p_1,p_2)$ into $B_V$ such that $\phi_p(x)> 0$ and $L_p \phi_p(x) = \Lambda(p) \phi_p(x)$ for all $p \in (p_1,p_2)$ and $x \in M$.
\end{theorem}

Assuming {\bf (G)} in Theorem \ref{thm kato1} can be satisfied with arbitrarily large positive values of $p_2$ and arbitrarily large negative values of $p_1$ then we obtain $\Lambda(p)$ defined and analytic for all $p \in \R$.    Since the value of $\Lambda(p)$ does not depend on the function $V$, it does not matter if larger and larger $|p_1|$ and $p_2$ require different functions $V$.   In this case the G\"{a}rtner-Ellis theorem, see for example \cite{DZ10}, implies $A_t/t$ satisfies the large deviations principle with speed $t$ and rate function $I(s) = \sup_{p \in \R} (ps-\Lambda(p))$.

Corollary \ref{cor Lam}(i) says that $\Lambda_t^x(p) \equiv \frac{1}{t} \log \E^x\left[e^{pA_t}\right]$ converges to $\Lambda(p)$ as $t \to \infty$ for each fixed $p$.  The following result gives a rate for this convergence and extends it to derivatives of the functions involved.  

\begin{theorem} \label{thm kato2}   Assume {\bf (H,P,G)} valid for two values $p_1 < p_2$ with same function $V$. 
For each $p_0 \in (p_1,p_2)$ and $x \in M$ there exist $\delta_1 > 0$, $t_0 < \infty$ and constants $C_0$, $C_1$ and $C_2$ such that for $k =0,1,2$
     \begin{equation}  \label{deriv}
     \left| (\Lambda_t^x)^{(k)}(p) - \Lambda^{(k)}(p)\right| \le \frac{C_k}{t} \quad \mbox{ for all } t \ge t_0 \mbox{ and } |p-p_0| < \delta_1.
    \end{equation} 
\end{theorem}

Proofs for Theorems \ref{thm kato1} and \ref{thm kato2} are given in Section \ref{sec vary proof}.

\subsection{Behavior for $p$ near 0} \label{sec near}

The large deviations of $A_t/t$ as $t \to \infty$ are determined, via the G\"{a}rtner-Ellis theorem, by the function $\Lambda(p)$ defined on the entire real line.  Here we concentrate on aspects of the limiting behavior of $A_t/t$ which are determined by the values of $\Lambda(p)$ for $p$ in an arbitrarily small neighborhood of $0$.  

Suppose first that the assumptions {\bf (H,P,G)} are satisfied at $p = 0$.  Putting $p = 0$ in Theorem \ref{thm Tt}(ii) implies $\{x_t: t \ge 0\}$ has the strong Feller property.  The eigenfunction $\phi_0$ can be chosen to be identically 1, and then 
$\widehat{\nu}_0$ is an invariant probability measure $\nu$, say, with a smooth density. We have $\int_M V d\nu < \infty$ and $\nu(U) > 0$ for every non-empty open $U \subset M$.  Since $\kappa_0 = 1$ the inequality \eqref{get} in Theorem \ref{thm KR} is a statement of $V$-exponential ergodicity for the Markov process $\{x_t:t \ge 0\}$.   The hypoellipticity condition {\bf (H)} implies that every invariant probability measure has a smooth density, and the positivity condition {\bf (P)} implies $\nu(U) > 0$ for all non-empty open $U \subset M$.  This implies $\nu$ is the unique invariant probability measure and that $\mbox{supp}(\nu) = M$.   It follows (see \cite[Prop 5.1]{IK}) that $P_t(x,\cdot)$ is absolutely continuous with respect to $\nu$ for all $t >0$ and all $x \in M$.

\s

The following Lemma shows that the growth condition {\bf (G)} valid for $p = \pm \delta$ is stronger than the integrability conditions typically made on the functions $Q$ and $q_1, \ldots, q_m$ in the additive functional $A_t$.  The proof is given in Section \ref{sec near proof}.    

\begin{lemma} \label{lem nu} Assume {\bf (H,P,G)} valid for two values $p = \pm \delta$ with same function $V$ for some $\delta > 0$.  Then $Q$, $YV/V$ and $R$ are all in $L^1(\nu)$. 
\end{lemma} 

In the setting of Lemma \ref{lem nu} write
    $$
    A_t = \int_0^tQ(x_s)ds + \sum_{j=1}^m \int_0^t q_j(x_s)dW_s^j = \int_0^tQ(x_s)ds + N_t
    $$
where $\{N_t: t \ge 0\}$ is a continuous local martingale with 
  $$
  \lim_{t \to \infty} \frac{1}{t}\langle N\rangle_t = \lim_{t \to \infty} \frac{1}{t} \int_0^t \left(\sum_{j=1}^r \bigl(q_j(x_s)\bigr)^2\right)ds  = \lim_{t \to \infty}  \frac{1}{t}\int_0^t R(x_s)ds.
  $$
Since $R \in L^1(\nu)$ and $P_t(x,\cdot)$ is absolutely continuous with respect to $\nu$ the almost sure ergodic theorem implies 
  $$ \lim_{t \to \infty} \frac{1}{t}\langle N\rangle_t  = \int_M R d\nu < \infty
  $$
$\PP^x$-almost surely.  It follows that $\frac{1}{t}N_t \to 0$ as $t \to \infty$ $\PP^x$-almost surely.  Applying the almost sure ergodic theorem again and using $Q \in L^1(\nu)$ gives
   \begin{equation} \label{lln}
   \lim_{t \to \infty} \frac{1}{t}A_t = \lim_{t \to \infty} \frac{1}{t}\int_0^t Q(x_s)ds = \int_M Q d \nu: = \lambda
   \end{equation}
$\PP^x$-almost surely for all $x \in M$.  In the setting of any of the motivating examples \eqref{v}, \eqref{v2} or \eqref{Tx} in Section \ref{sec intro} with $A_t = \log\|v_t\| - \log \|v_0\|$, this can be recognized as the Furstenberg-Khas'minskii formula for the top Lyapunov exponent $\lambda$, see \cite{Arn98}.

\s

In the following result note that the smoothness and convexity of $\Lambda(p)$ for $-\delta <p< \delta$ implies the existence of $\Lambda'(0) \in \R$ and $\Lambda''(0) \ge 0$.

\begin{theorem} \label{thm clt} Assume {\bf (H,P,G)} valid for two values $p = \pm \delta$ with same function $V$ for some $\delta > 0$.  Then
   
   (i)  $\Lambda(0) = 0$ and $\Lambda'(0) = \lambda $.
   
   (ii) For all $x \in M$ the distribution of $$
     \sqrt{t}\left(\frac{A_t}{t}-\Lambda'(0) \right)
     $$
under $\PP^x$ converges as $t \to \infty$ to the normal distribution $N(0,\Lambda''(0))$ with mean 0 and variance $\Lambda''(0)$.
     
\end{theorem}

In Section \ref{sec near proof} we present two methods of proof for Theorem \ref{thm clt}.  The first method, in Section \ref{sec near abstract}, uses the more abstract properties $\Lambda(p)$ in Theorem \ref{thm kato2} and is based very closely on results of Liming Wu \cite{WuL95}.  It gives rise also to the following moderate deviations result, see \cite[Theorem 1.2]{WuL95}.

\begin{theorem} \label{thm mod dev}  Assume {\bf (H,P,G)} valid for two values $p = \pm \delta$ with same function $V$ for some $\delta > 0$.   Assume $\Lambda''(0) > 0$. Suppose $b_t/\sqrt{t} \to \infty$ and $b_t/t \to 0$ as $t \to \infty$.  Then $(A_t-t \Lambda'(0))/b_t$ satisfies the large deviation principle with speed $b_t^2/t$ and rate function $I(s) = s^2/(2 \Lambda''(0))$.

\end{theorem}

For example in Theorem \ref{thm mod dev} we could take $b_t = t^\beta$ for $1/2 < \beta < 1$.

\s

The second method of proof of Theorem \ref{thm clt} uses properties of the eigenfunctions $\phi_p$ associated with $\Lambda(p)$.  More specifically it uses the real analyticity of $p \to \Lambda(p)$ and $p \to \phi_p$ in the equation
  \begin{equation} \nonumber
   \Bigl(L+ p \sum_{j=1}^m q_j X_j + pQ + \frac{1}{2}p^2 \sum_{j=1}^m q_j^2\Bigr) \phi_p(x) = \Lambda(p)\phi_p(x).
  \end{equation}
In particular it gives the existence of the function $\Phi_1 = \dfrac{\partial \phi_p}{\partial p}\big|_{p=0}:M \to \R$ satisfying the Poisson equation 
  \begin{equation} \nonumber
 Q(x) + L\Phi_1(x) = \Lambda'(0)
 \end{equation} 
with the property that $(\Phi_1)^q \in B_V \subset L^1(\nu)$ for all $q \ge 1$. 
Details are given in Section \ref{sec near eigen}.  The use of eigenfunctions also allows the following characterization of the degenerate case when $\Lambda''(0) = 0$.

\begin{theorem} \label{thm linear} Assume {\bf (H,P,G)} valid for two values $p = \pm \delta$ with same function $V$ for some $\delta > 0$.  The following are equivalent

(i) $\Lambda''(0) = 0$.

(ii) There exists $\Phi_1 \in B_V$ such that 
  $X_j \Phi_1(x) + q_j(x)  = 0$ for all $j \ge 1$ and $L\Phi_1(x) + Q(x)  = \Lambda'(0)$.

(iii)  There exists $\Phi_1 \in B_V$ such that $A_t = t \Lambda'(0) + \Phi_1(x_0) - \Phi_1(x_t)$.

(iv) $\Lambda(p) = p \Lambda'(0)$ for all $|p| <\delta$.
\end{theorem}

\section{Examples} \label{sec examples}  We consider briefly two examples in which the results above can be applied to the linearization of a nonlinear stochastic differential equation.

\subsection{Langevin example}   The system 
    \begin{align*}
     dx_t & = y_t dt\\
     dy_t & = (ax_t-bx_t^3 - \beta y_t)dt+ \sigma dW_t
     \end{align*}
with parameters $a \in \R$, $b > 0$, $\beta > 0$ and $\sigma \neq 0$ is a particular case of the Langevin family studied by Mattingly, Stuart and Higham \cite{MSH02} and Ferr\`{e} and Stoltz \cite{FS20}.  Whereas those papers consider ergodic and large deviation behavior of the process $\{(x_t,y_t): t \ge 0\}$ in $\R^2$, here we consider also the linearization
   $$
   dv_t = \begin{bmatrix} 0 & 1 \\ a-3bx_t^2 & -\beta \end{bmatrix}v_tdt.  
   $$
Putting $v_t = \|v_t\| \begin{bmatrix} \cos \theta_t \\ \sin \theta_t \end{bmatrix}$ gives  
   \begin{align*}
    d \theta_t & = \bigl(-\sin^2 \theta_t + (a-3bx_t^2) \cos^2 \theta_t - \beta \sin \theta_t \cos \theta_t \bigr)dt \\
    d \log\|v_t\| &  = \bigl((1+a-3bx_t^2)\sin \theta_t \cos \theta_t - \beta \sin^2 \theta_t\bigr)dt.
   \end{align*}
This gives a diffusion $\{(x_t,y_t,\theta_t); t \ge 0\}$ on  $\R^2 \times \mathbb{S}^1$, and then $A_t = \log\|v_t\|-\log \|v_0\|$ has $Q(x,y,\theta) = (1+a-3bx^2)\sin \theta \cos \theta - \beta \sin^2 \theta$ and $q_j \equiv 0$. 

Let $\widehat{L}$ denote the generator of $\{(x_t,y_t): t \ge 0\}$ and let $L$ denote the generator of $\{(x_t,y_t,\theta_t): t \ge 0\}$.  Notice that $|Q(x,y,\theta)| \le A+Bx^2$ for constants $A$ and $B$, and define 
 \begin{equation} \label{V}
  V(x,y) = \exp\left\{ \gamma\left(- \frac{ax^2}{2} + \frac{bx^4}{4} + \frac{y^2}{2}\right) + \e xy \right\}.  
  \end{equation}
  It is shown in \cite{MSH02} that there exist $\gamma > 0$ and $\e > 0$ such that $\widehat{L}V(x,y)/V(x,y) \to - \infty$ as $\|(x,y)\| \to \infty$.  Since $V$ does not depend on $\theta$, essentially the same calculations can be used to show that 
     $$
    \frac{LV(x,y)}{V(x,y)} + pQ(x,y,\theta) \le  \frac{\widehat{L}V(x,y)}{V(x,y)} +p(Ax^2+B) \to -\infty
   $$   as  
   $\|(x,y)\| \to \infty$      
for all $p \in \R$, $0 < \gamma < 2\beta/\sigma^2$ and sufficiently small $\e > 0$. It follows that the growth condition {\bf (G)} is satisfied with $V$ as in \eqref{V} for all $p \in \R$.  

The conditions {\bf (H)} and {\bf (P)} for the process $\{(x_t,y_t: t \ge 0\}$ are verified in \cite{MSH02}.  It is a non-trivial matter to extend these conditions to the 3-dimensional system $\{(x_t,y_t,\theta_t): t \ge 0\}$, and this will be done elsewhere.

\subsection{Nonlinear autoparametric system}  In Baxendale and Sri Namachchivaya \cite{BSri} the stability of the single mode solution of a nonlinear autoparametric system is given by the limiting behavior of $\frac{1}{t}\log \|u_t\|$ as $t \to \infty$ for the 4-dimensional system
   \begin{align*}
   dv_t & = \begin{bmatrix} 0 & 1 \\-\chi^2 & -2 \zeta_1 \end{bmatrix}dv_t + \begin{bmatrix}0 \\ \sigma \end{bmatrix}dW_t \\
   du_t & = \begin{bmatrix} 0 & 1 \\ -\kappa^2 -\chi^2 v_t^{(1)} - 2 \zeta_1 v_t^{(2)} &  -2 \zeta_2 \end{bmatrix}du_t + \begin{bmatrix}0 & 0 \\ \sigma & 0\end{bmatrix}u_t dW_t
    \end{align*}
(see \cite[Equations 2.6, 2.7, 2.8]{BSri}; here $v_t$ plays the role of $x_t \in \R^2$ and $u_t$ plays the role of $v_t$).   Putting $u_t = \|u_t\| \begin{bmatrix} \cos \psi_t \\ \sin \psi_t \end{bmatrix}$ gives  
   \begin{align*}
  d\psi_t & = h(v_t,\psi_t)dt+ \sigma \cos^2 \psi_t dW_t \\
  d \log \|u_t\| &= Q(v_t,\psi_t)dt+ \sigma \sin \psi_t \cos \psi_t dW_t
  \end{align*}
where 
  \begin{align*}
  h(v,\psi) & = -1+\bigl(1-\kappa^2-\chi^2 v^{(1)}-2\zeta_1 v^{(2)}\bigr)\cos^2 \psi - 2\zeta_2 \sin \psi \cos \psi -\sigma^2 \sin \psi \cos^3 \psi \\
  Q(v,\psi) & = \bigl(1-\kappa^2-\chi^2 v^{(1)}-2\zeta_1 v^{(2)}\bigr)\sin \psi \cos \psi - 2\zeta_2 \sin^2 \psi+ \frac{\sigma^2}{2} \cos^2 \psi \cos 2 \psi.
  \end{align*}     
This gives a diffusion $\{(v_t,\psi_t); t \ge 0\}$ on  $\R^2 \times \mathbb{S}^1$, and then $A_t = \log\|u_t\|$ has $Q(v,\psi)$ as above and $q(v,\psi) = \sigma \sin \psi \cos \psi$.  

The conditions {\bf (H)} and {\bf (P)} are verified in \cite[Section 8]{BSri}.  
 Here we will verify {\bf (G)}.  Let $\widehat{L}$ denote the generator of $\{v_t: t \ge 0\}$ and let $L$ denote the generator of $\{(v_t,\psi_t): t \ge 0\}$.  Write $A =  \begin{bmatrix} 0 & 1 \\-\chi^2 & -2 \zeta_1 \end{bmatrix}$ and $\overline{\sigma} = \begin{bmatrix} 0 \\ \sigma \end{bmatrix}$.  
  Since $A$ has eigenvalues with strictly negative real parts, there exists an invertible real $2 \times 2$ matrix $C$ and $\delta >0$ such that $\langle CAv,Cv \rangle \le -\delta \|Cv\|^2$ for all $v \in \R^2$.  Choose $U(v) = \|Cv\|^2$ and $V(v) = e^{\gamma U(v)}$ for some $\gamma > 0$.   
   If $X(v,\psi)$ denotes the vector field multiplying $dW_t$ in the SDE for $\{(v_t,\psi_t): t \ge 0\}$ then $XU(v,\psi) = 2 \langle Cv,C \overline{\sigma}\rangle$.   Then for each $p \in \R$, using \eqref{LpV} we have
  \begin{align*}
    \frac{(L_pV)(v,\psi)}{V(v)} & = 
    \gamma\widehat{L}U(v) + \frac{1}{2}\gamma^2 \bigl(XU(v)\bigr)^2 + p \gamma  q(v,\psi)XU(v) +pQ(v,\psi) + \frac{1}{2}p^2 \bigl(q(v,\psi)\bigr)^2 \\
        & = 
    \gamma\Bigl(2 \langle CAv,Cv\rangle + \langle C \overline{\sigma}, C\overline{\sigma} \rangle \Bigr) 
        + 2 \gamma^2 \langle Cv,C \overline{\sigma}\rangle^2 + 2 \gamma p q(v,\psi)\langle Cv,C \overline{\sigma}\rangle \\
    & \quad + pQ(v,\psi)+ \frac{1}{2}p^2 \bigl(q(v,\psi)\bigr)^2 \\
        & \le -2\gamma \delta \|Cv\|^2 + 2 \gamma^2\|C\overline{\sigma}\|^2 \|Cv\|^2 + {\cal O}(\|v\|) \\
    & \to -\infty     \end{align*}     
as $\|v\| \to \infty$ so long as $\gamma < \delta/\|C\overline{\sigma}\|^2$.
  Therefore {\bf (G)} is satisfied for all $p \in \R$ with $V(v) = e^{\gamma \|Cv\|^2}$ for any $\gamma$ satisfying $0 <\gamma < \delta/\|C\overline{\sigma}\|^2$.

\section{Some examples of asymptotic behavior of $\Lambda(p)$.}   \label{sec asymp}
In the case that $M$ is compact and $q_j \equiv 0$ for all $j \ge 1$ it is clear that
   $$
  \frac{\Lambda(p)}{p} \le \max\{Q(x): x \in M\} < \infty
   $$
for all $p$.  Moreover, when $M$ is compact and $L$ is strictly elliptic then
    $$
    \lim_{p \to \pm \infty} \frac{\Lambda(p)}{p^2} = \frac{1}{2} \inf_{g \in C^\infty(M)} \sup_{x\in M} \sum_{j=1}^m \Bigl(X_jg(x)-q_j(x)\Bigr)^2 < \infty, 
    $$
see \cite[Theorem(1.25)]{BS88}.  Not surprisingly, the asymptotic behavior of $\Lambda(p)$ as $p \to \infty$ can be very different in the non-compact case.  In this section we give some particular results when the underlying diffusion process $\{x_t: t \ge 0\}$ is the pitchfork bifurcation with additive noise
   \begin{equation} \label{pitch}
    dx_t= (ax_t-bx_t^3)dt+ \sigma dW_t.
    \end{equation} 
There is no linearization carried out in these examples.  Instead we make arbitrary choices for the functions $Q(x)$ and $q(x)$.

\begin{proposition} \label{prop asymp} For the system \eqref{pitch} with $b >0$ and $\sigma > 0$

(i) if $\displaystyle{A_t = \int_0^t x_s^2\,ds}$ then $
     \lim_{p \to \infty} \dfrac{\Lambda(p)}{p^{3/2}} = \left(\dfrac{2}{3}\right)^{3/2} \dfrac{\sigma}{b};
     $
     
     \s
     
(ii) if $\displaystyle{A_t = \int_0^t x_s^4\,ds}$ then $
     \lim_{p \to \infty} \dfrac{\Lambda(p)}{p^3} = \dfrac{16}{27}\left(\dfrac{\sigma}{b}\right)^4;
     $
 
 \s  
 
(iii) if  $\displaystyle{A_t = \int_0^t x_s\,dW_s}$  then $
     \lim_{p \to \infty} \dfrac{\Lambda(p)}{p^3} = \dfrac{16 \sigma}{27b};$

\s

(iv) if $\displaystyle{A_t = -\int_0^t x_s\,dW_s}$ then 
    $ 
  \displaystyle{ \frac{(\max(a,0))^2}{4b\sigma} - \frac{\sigma}{2} \le  \lim_{p \to \infty} \frac{\Lambda(p)}{p}  \le \frac{(\max(a,0))^2}{4b\sigma} + \frac{\sigma}{2}}.
   $
   
 \end{proposition}      

\begin{remark} \label{rem corr} {\rm It can be seen in (iii) and (iv) that the relationship between the noise in \eqref{pitch} and the noise in $A_t$ is important.  The condition that $\sigma > 0$, rather than just $\sigma \neq 0$, is important in distinguishing cases (iii) and (iv).   Suppose more generally 
$A_t = \int_0^t x_s d\widetilde{W}_s$ where $\{W_t: t \ge 0\}$ and $\{\widetilde{W}_t: t \ge 0\}$ are standard Wiener processes with covariation $\langle W, \widetilde{W} \rangle_t = \rho t$.  The correlation coefficient $\rho$ occurs frequently in stochastic volatility models in mathematical finance, see \cite{heston}, \cite{fouque}.  This fits within the general setting of an SDE \eqref{x} and an additive functional \eqref{A} by writing
    \begin{align*}
    dx_t & = (ax_t-bx_t^3)dt+ \sigma dW_t^1\\
    A_t& = \rho \int_0^t x_sdW_s^1+ \sqrt{1-\rho^2}\int_0^t x_s dW_t^2
    \end{align*}  
where $\{W_t^1:t \ge 0\}$ and $\{W_t^2: t \ge 0\}$ are independent standard Wiener processes.  Then $Y(x) = \sum_{j=1}^2 q_j(x)X_j(x) = \rho \sigma x \frac{d}{dx}$ and $R(x) = \sum_{j=1}^2 \bigl(q_j(x)\bigr)^2 = x^2$.  The extreme cases $\rho = \pm 1$ appear as cases (iii) and (iv) above.  For $-1 < \rho \le 1$ we have
  $$
   \lim_{p \to \infty} \frac{\Lambda(p)}{p^3} = \frac{\Bigl(\sqrt{3+\rho^2}+ 2 \rho\Bigr)^2\Bigl(\sqrt{3+\rho^2}-\rho \Bigr)\sigma}{27b}.
  $$  
} \end{remark}

\s

The proofs of these limits appear in Section \ref{sec asymp proof}.  The method of proof relies on two observations which are valid for any system satisfying {\bf (H,P,G)} for some fixed $p$ with a given $V$.
 \begin{itemize}
   \item  Suppose $\sup \left\{\frac{L_pV(x)}{V(x)}: x \in M\right\} = \Gamma$.  Lemma \ref{lem Vest} gives $\|T_t^p\|_V \le e^{t \Gamma}$ for all $t \ge 0$ and so $\Lambda(p)= \log \kappa_p \le \Gamma$.

\item Suppose $g \in B_V \cap C^2(M)$ is non-negative, not identically 0, and satisfies $L_p g \ge \widetilde{\Lambda} g$.  By Remark \ref{rem pos} we have $\Lambda(p) \ge \widetilde{\Lambda}$. 
\end{itemize}
These two observations can be used to obtain upper and lower bounds on $\Lambda(p)$.  For the systems in Proposition \ref{prop asymp} the asymptotic behaviors of the upper and lower bounds happen to coincide.  
  
\section{Proofs for Section \ref{sec fix}} \label{sec fix proof}

Within this Section only we break down the growth assumption {\bf (G)} into its invidual parts {\bf (0)}, \ldots ,{\bf (3)} to see where each part of {\bf (G)} is used.


\subsection{Non-explosion and bounds on $T_t^p$}    \label{sec growth}

 For $k \ge 1$ define $D_k = \{x \in M: V(x) \le k\}$.  Then $D_k$ is compact and $D_k \subset \mbox{int}(D_{k+1})$ and $D_k \nearrow M$. 
  For each initial point $x \in M$ the process $\{x_t: t \ge 0\}$ is well defined up to the stopping time $\tau_k = \inf\{t \ge 0: x_t \not\in D_k\}$.  The first result is due to Meyn and Tweedie \cite[Theorem 2.1]{MTIII}.
 
\begin{lemma}  \label{lem stop} Assume {\bf (0)}.  
There exists $\Gamma_0 < \infty$ such that 
     $$
     \PP^x(\tau_k \le t) \le k^{-1} V(x) e^{\Gamma_0 t}
     $$
whenever $V(x) < k$.    In particular $\PP^x(\tau_k \to \infty \mbox{ as } k \to \infty) = 1$ for all $x$, so that that the process $\{x_t: t \ge 0\}$ is non-explosive, that is, the solution exists for all time $t \ge 0$.
 
\end{lemma}           

\n{\bf Proof.}  Define $\Gamma_0 = \sup_{x \in M}LV(x)/V(x) \in [0,\infty)$.  Fix $k \ge 1$ and define $M_t = V(x_{t \wedge \tau_k})e^{- \Gamma_0 (t \wedge \tau_k)}$.  Adapting the method of proof of Lemma \ref{lem Vest} to the condition {\bf (0)} we deduce that $M_t$ is a supermartingale.  
For $V(x) < k$, on the set $\tau_k \le t$ we have $M_t = V(x_{\tau_k})e^{ -\Gamma_0 \tau_k} = k e^{ -\Gamma_0 \tau_k} \ge k e^{-\Gamma_0 t}$.    
Therefore 
   $$k e^{-\Gamma_0 t}\PP^x(\tau_k \le t) \le  \E^x\left[M_t \mathbf{1}_{\tau_k \le t} \right] \le  \E^x\left[M_t\right] \le V(x),
   $$     
and the result follows directly.  \qed

\s

The next result uses only a weakened version of {\bf (1)} in which $L_pV/V$ is bounded above.  

\begin{lemma}  \label{lem Vest}  Assume {\bf (0,1)}.  There exists $\Gamma$ such that 
   \begin{equation} \label{Vest} 
      T_t^p V(x) \le e^{\Gamma t} V(x)
      \end{equation}
for all $t \ge 0$ and $x \in M$.  Further, $T_t^p$ is a well-defined bounded linear operator on $B_V$ with operator norm $\|T_t^p\|_V \le e^{\Gamma t}$.
  \end{lemma}

\n{\bf Proof.}  The proof follows closely Wu \cite[Corollary 2.2]{WuL01} and Ferr\'{e} and Stoltz \cite[Lemma 6.2]{FS20}.  The condition {\bf (1)} implies 
   \begin{equation} \label{Gamma}
   \Gamma: = \sup_{x \in M}\left(\frac{L_pV(x)}{V(x)} \right) = \sup_{x \in M} \left(\frac{LV(x)+pYV(x)}{V(x)} +pQ(x) + \frac{p^2}{2}R(x)\right)
   \end{equation} 
is finite.  Let $M_t =  V(x_t)e^{p A_t}e^{-\Gamma t}$, then It\^{o}'s formula gives
\begin{align}
    dM_t & = \left( LV(x_t) + pYV(x_t) +pQ(x_t)V(x_t)+\frac{p^2}{2}R(x_t)V(x_t) - \Gamma V(x_t) \right)e^{pA_t} e^{-\Gamma t} dt  \nonumber \\
   & \qquad        + \sum_{j=1}^m \bigl(X_jV(x_t)+pq_j(x_t)V(x_t)\bigr)e^{p A_t}e^{-\Gamma t}dW_j(t) \nonumber \\  
     & = \left( \frac{(L_pV(x_t)}{V(x_t)}  - \Gamma \right)M_t dt  \nonumber\\
      & \qquad        + \sum_{j=1}^m \bigl(X_jV(x_t)+p q_j(x_t)V(x_t)\bigr)e^{p A_t}e^{-\Gamma t}dW_j(t) \label{dMt}
       \end{align}
so that $M_t$    
is a non-negative local supermartingale and hence a supermartingale.  Therefore
  $$
  T_t^pV(x) =  e^{\Gamma t} \E^x[M_t] \le e^{\Gamma t}\E^x[M_0] = e^{\Gamma t}V(x)
  $$
and \eqref{Vest} is proved.  For $g \in B_V$ we have $|g(x_t)| \le \|g\|_V V(x_t)$ and so 
  \begin{align*}
      \E^x\left|g(x_t)  e^{p A_t}\right| 
       \le \|g\|_V \E^x\left[V(x_t)e^{ pA_t}\right]
       = \|g\|_V T_t^pV(x)  \le \|g\|_V e^{\Gamma t}V(x).
    \end{align*}
Therefore $T_t^pg(x)$ is well defined and $|T_t^pg(x)| \le \|g\|_V e^{\Gamma t}V(x)$.  
\qed

\s

The following result, stated without proof, uses a similar argument with the other growth conditions. 
 
\begin{lemma}  \label{lem Vest234} 
(i) Assume {\bf (0,2)}.   There exists $\beta > 1$ and $\Gamma_2 \ge 0$ such that 
   \begin{equation} \label{Vest2} 
     \E^x\left[V(x_{t \wedge \tau}) e^{\beta p A_{t \wedge \tau}}\right] \le e^{\Gamma_2 t} V(x)
      \end{equation}
for all $t \ge 0$, $x \in M$ and stopping times $\tau$.

 (ii) Assume {\bf (0,3)}.  There exists $\beta > 1$ and $\Gamma_3 \ge 0$ such that 
   \begin{equation} \label{Vest3} 
      \E^x\left[V^\beta(x_{t \wedge \tau}) e^{p A_{t \wedge \tau}}\right]\le e^{\Gamma_3 t} V^\beta(x)
      \end{equation}
for all $t \ge 0$, $x \in M$ and stopping times $\tau$.
  \end{lemma}

 \s
 
The next result is the only place where the full strength of $L_pV(x)/V(x) \to -\infty$ as $x \to\infty$ in {\bf (1)} is used.  It will be important in the proof that $T_t^p$ is a compact operator, see Lemma \ref{lem comp}.  A similar result with $q_j \equiv 0$ appears as part of \cite[Lemma 6.2]{FS20}.

  \begin{lemma} \label{lem Vest12}  Assume {\bf (0,1,2)}.   For each $t >0$ we have 
       \begin{equation} \label{TtVV}
        \frac{T_t^pV(x)}{V(x)} \to 0 \quad \mbox{ as } x \to \infty.
        \end{equation}
\end{lemma} 

\n{\bf Proof.} Assumption {\bf (1)} implies $\dfrac{L_p V(x)}{V(x)} \to -\infty$ as $x \to \infty$, so there exist sequences $a_n \to \infty$ and $b_n \in \R$  such that $a_nV(x) + L_pV(x) \le b_n$ for all $x \in M$.    Let $M_t =  e^{a_n t}V(x_t)e^{p A_t}$, then (adapting the calculation in \eqref{dMt})
  \begin{align*}
    dM_t & = \bigl(a_nV(x_t) +L_p V(x_t) \bigr)e^{a_n t}e^{p A_t} dt\\
    & \quad +  \sum_{j=1}^m \bigl(X_jV(x_t)+ p q_j(x_t)V(x_t)\bigr)e^{a_n t}e^{p A_t} dW_j(t).
       \end{align*}
The integrals of the noise terms are local martingales, so with the stopping times $\tau_k \to \infty$ from Lemma \ref{lem stop} we have 
     \begin{align*}
     \E^x[M_{t \wedge \tau_k}] 
     & \le  V(x)+ \E^x\left[\int_0^{t \wedge \tau_k}\bigl(a_nV(x_s) + L_p V(x_s) \bigr)e^{a_n s}e^{p A_s}ds\right] \\
      & \le   V(x)+b_n \E^x\left[\int_0^{t\wedge \tau_k} e^{a_n s} e^{p A_s}ds\right] \\
          & \le  V(x)+ b_n \int_0^t e^{a_n s} (T_s^p \mathbf{1})(x)ds. 
          \end{align*}
Letting $\tau_k \to \infty$ and using Fatou's lemma gives
   $$
   e^{a_n t} (T_t^p V)(x)  = \E^x[M_t]  
     \le  V(x)+ b_n \int_0^t e^{a_n s} (T_s \mathbf{1})(x) ds.
     $$
By Lemma \ref{lem Vest234}(i), using assumption {\bf (2)}, there exist $\beta > 1$ and $\Gamma_2 < \infty$ such that
  $$
    (T_t^p \mathbf{1})(x) = \E^x\left[e^{p A_t}\right] \le \left(\E^x\left[e^{ \beta p A_t}\right]\right) ^{1/\beta} \le  \left(\E^x\left[V(x_t)e^{ \beta p A_t}\right]\right) ^{1/\beta} \le \left(e^{\Gamma_2 t} V(x)\right)^{1/\beta}
    .$$ 
Therefore 
    $$
   e^{a_n t} (T_t^p V)(x)  \le  V(x)+b_n \int_0^t e^{(a_n+\Gamma_2/\beta) s} (V(x))^{1/\beta} ds \le V(x)+ \frac{b_n e^{(a_n+\Gamma_2/\beta)t}}{a_n+\Gamma_2/\beta}(V(x))^{1/\beta}.
     $$
Since $(V(x))^{1/\beta}/V(x) = V(x)^{-(\beta-1)/\beta} \to 0$ as $x \to \infty$ we get
     $$
     \limsup_{x \to \infty} \frac{T_t^pV(x)}{V(x)} \le e^{-a_n t}
     $$
and result follows by letting $n \to \infty$. \qed

\s
\begin{remark} In this result the assumption {\bf (2)} can be replaced by the assumption that $L_p \widetilde{V}(x)/\widetilde{V}(x)$ is bounded  above by $\widetilde{\Gamma}$ for some $\widetilde{V} \ge 1$ with $\widetilde{V}(x)/V(x) \to 0$ as $x \to \infty$.  See \cite[Assumption 4]{FRS21}.   For then $T_t^p \mathbf{1}(x) \le e^{\widetilde{\Gamma}t }\widetilde{V}(x)$ and the proof can be completed as above
\end{remark} 

\begin{remark}  It is easy to verify the semigroup property $T_t^p \circ T_s^p = T_{t+s}^p$ for $t,s > 0$.  However the result of Lemma \ref{lem Vest12} implies that $\|T_t^pV - V\|_V \ge 1$ for all $t >0$, so that $\{T_t^p: t \ge 0\}$ is not a strongly continuous semigroup of operators.   It does not have a generator with a dense domain in $B_V$ in the usual sense.
\end{remark}

\subsection{Regularity of the semigroup}  

In this section we prove parts (ii) and (iii) of Theorem \ref{thm Tt}.   
Here is a brief summary of the sequence of results.
\begin{itemize}
   \item Use smoothness of the $X_j$ and $Q$ and $q_j$ and a truncation argument and a Girsanov transformation to show 
       $$
       (L_{p,k}-\frac{\partial}{\partial t})(T_{t,k}^p g)(x) = 0, \quad \mbox{ for }g \in C^\infty_c(M);
       $$
this is Lemma \ref{lem FK k}.  The important idea is that a smooth stochastic flow preserves the smoothness of the initial condition, see Ikeda and Watanabe \cite[Chapter V.]{IWbook}.

\item Let $k \to\infty$ in the sequence of truncations to show
        $$
       (L_p-\frac{\partial}{\partial t})(T_t^p g)(x) = 0
       $$
in the sense of distributions for $g \in C^\infty_c(M)$; this is Proposition \ref{prop FK phi}

\item  Adapt a result of Ichihara and Kunita \cite{IK}, using {\bf (H)} to get
      $$
      T_tg(x) = \int \rho_t(x,y)g(y)dy
      $$
for $g \in B_V$, where $\rho:(0,\infty) \times M \times M \to \R$ is $C^\infty$ and $dy$ denotes Riemannian measure on $M$; this is Proposition \ref{prop kernel}.

\item Use continuity of $(x,y) \to \rho_t(x,y)$ and $T_tV(x)/V(x) \to 0$ as $x \to \infty$ (Lemma \ref{lem Vest12}) to show $T_t$ is compact for $t >0$, and that $x \mapsto T_tg(x)$ is continuous for all $g \in B_V$ and $t > 0$; this is Lemma \ref{lem comp} and it proves Theorem \ref{thm Tt}(iii).

\item  Approximate continuous $g \in B_V$ by $\phi_n \in C^\infty_c(M)$ to extend $(L_p-\frac{\partial}{\partial t})(T_t^pg)(x) = 0$ from $g \in C^\infty_c(M)$ to $g \in B_V$.   This is Proposition \ref{prop FK g} and it proves Theorem \ref{thm Tt}(ii)

\end{itemize}

Recall the sets $D_k = \{x \in M: V(x) \le k\}$ from Section \ref{sec growth}.  For each $k$ choose a $C^\infty$ function $\chi_k:M \to [0,1]$ such that $\chi_k(x) = 1$ for $x \in D_k$ and $\mbox{supp}(\chi_k) \subset \mbox{int}(D_{k+1})$.  
Define the truncated vector fields $X_{j,k}(x) = \chi_k(x) X_j(x)$ for $1 \le j \le m$ and $X_{0,k}(x) = \chi_k^2(x)X_0(x) - \frac{1}{2}\sum_{j=1}^m \chi_k(x)(X_j\chi_k)(x)X_j(x)$.  The resulting truncated SDE
given by 
   \begin{equation} \label{xk}
    dx_{t,k} = X_{0,k}(x_{t,k})dt + \sum_{j=1}^m X_{j,k}(x_{t,k}) \circ dW_j(t)
  \end{equation} 
has a unique solution $\{x_{t,k}: t  \ge 0\}$ valid for all $t \ge 0$.  It has generator $L^{(k)}$, say, given by $L^{(k)}g(x) = \chi_k^2(x)Lg(x)$.    Define the truncated functions $Q_k(x) = \chi_k^2(x)Q(x)$ and $q_{j,k}(x) = \chi_k(x)q_j(x)$, and then the truncated additive functional 
    \begin{equation} \label{Ak}
    A_{t,k} = \int_0^t Q_k(x_{s,k})ds+ \sum_{j=1}^m \int_0^t q_{j,k}(x_{s,k}) dW_j(s),
  \end{equation}
and the truncated operator   
     $$
     T_{t,k}^p g(z) = \E^x\left[g(x_{t,k})e^{p A_{t,k}}\right].
     $$
     
\begin{lemma} \label{lem FK k}  Suppose $\phi \in C^\infty_c(M)$, the space of $C^\infty$ functions on $M$ with compact support, and define $u_k(t,x) = (T^p_{t,k}\phi)(x)$.  Then $u_k \in C^\infty((0,\infty) \times M)$ and
  \begin{equation} \label{FK k}
  \left(\chi_k^2(L+pY+pQ+\frac{p^2}{2}R) - \frac{\partial }{\partial t}\right)u_k(t,x)= 0
  \end{equation}
for all $t > 0$ and $x \in M$.  
\end{lemma}

\n{\bf Proof.}   Define 
 $$
   N_{t,k}  = \sum_{j=1}^m \int_0^t p q_{j,k}(x_{s,k})dW_j(s) - \frac{1}{2}\sum_{j=1}^m \int_0^t p^2 \bigl(q_{j,k}(x_{s,k})\bigr)^2ds 
   $$
Since the functions $q_{j,k}$ are bounded, then $e^{N_{t,k}}$ is a true martingale and we may define a probability measure $\PP_k^x$ on the canonical space of paths in $M$ by 
   $$
        \frac{d \PP_k^x}{d\PP^x} \Big|{{\cal F}_t} = e^{N_{t,k}}.
        $$   
Then  
   \begin{align}
      u_k(t,x)  & = \E^x \left[\phi(x_{t,k}) \exp\left\{\int_0^t p Q_k(x_{s,k})ds + \sum_{j=1}^m \int_0^t p q_{j,k}(x_{s,k})dW_j(s)\right\} \right] \nonumber  \\
      & = \E^x \left[\phi(x_{t,k}) \exp\left\{\int_0^t \left( p Q_k(x_{s,k})+ \frac{p^2}{2}\sum_{j=1}^m (q_{j,k}(x_{s,k}))^2\right)ds \right\} e^{N_{t,k}}\right] \nonumber \\
      & =\E_k^x  \left[\phi(x_{t,k}) \exp\left\{\int_0^t \chi_k^2(x_{s,k})\left(p Q(x_{s,k})+ \frac{p^2}{2}R(x_{s,k})\right)ds \right\}\right] 
      \label{uk}
      \end{align}
where $\E_k^x$ denotes expectation with respect to $\PP_k^x$.  By Girsanov's theorem the law of the diffusion $\{x_{t,k}: t \ge 0\}$ under $\PP_k^x$ is generated by an SDE of the form \eqref{xk} with drift vector field $X_{0,k}(x)$ replaced by $X_{0,k}(x) + \sum_{j=1}^m  p q_{j,k}(x)X_{j,k}(x) = X_{0,k} +  p\chi_k^2(x)Y(x)$.  It has generator $L_k + p\chi_k^2 Y = \chi_k^2(L+pY)$.  
Since the vector fields $X_{0,k} + \chi_k^2 Y$ and $X_{j,k}$ for $1 \le j \le m$ are all $C^\infty$ with compact support, the stochastic flow generated by this SDE consists of $C^\infty$ diffeomorphisms of $M$, see Ikeda and Watanabe \cite[Sect V.2]{IWbook}.   Since the function $\chi_k^2(x)(Q(x)+\frac{1}{2}R(x))$ also has compact support, then the result follows as an application of \cite[Theorem V.3.2]{IWbook} to the expression \eqref{uk}. \qed

\begin{remark}  
In this paper we use $L$, and then $L+pY$ and $\chi_k^2(L+pY)$, for the action of the generator of the relevant Markov process on $C^2$ functions.  The Feynman-Kac formula and Girsanov's theorem together give a result similar to \eqref{FK k} in which $\chi_k^2(L+pY)$ acts as the extended generator and there is no guarantee of smoothness of $u_k(t,x)$.
\end{remark}

Next we let $k \to \infty$ in the truncation.  

\begin{lemma} \label{lem kinf} Assume {\bf (0,1,2)}.  Suppose $g \in B_V$ is bounded.  Then $T_{t,k}^p g(x) \to T_t^p g(x)$ uniformly on compact subsets of $[0,\infty) \times M$ as $k \to \infty$.
    \end{lemma}

\n{\bf Proof.}   We will show uniform convergence on a set of the form $[0,T] \times D_{k_0}$ for some fixed $T$ and $k_0$.  Without loss of generality we may assume that $|g(x)| \le 1$ for all $x \in M$.  Recall the stopping times $\tau_k = \inf\{t \ge 0 : x_t \not\in D_k\}$.  
Since $\{x_{s,k}: 0 \le s \le t\} = \{x_s: 0 \le s \le t\}$ and $\{A_{s,k}: 0 \le s \le t\} = \{A_s: 0 \le s \le t\}$ on the set $\{\tau_k > t\}$ we have 
 \begin{align*}
  \big|T^p_{t,k}g(x) - T^p_t g(x)\big|
   & = \big| \E^x\bigl[g(x_{t,k}) e^{p A_{t,k}}\bigr] -  \E^x\bigl[g(x_t) e^{p A_t}\bigr] \big| \\
    & = \big| \E^x\bigl[g(x_{t,k}) e^{p A_{t,k}}1_{\{\tau_k\le t\}}\bigr] -  \E^x\bigl[g(x_t) e^{p A_t}1_{\{\tau_k\le t\}}\bigr] \big| \\
    & \le  \E^x\left[ e^{p A_{t,k}}1_{\{\tau_k\le t\}} \right] +  \E^x\left[ e^{p A_t}1_{\{\tau_k \le t\}}\right]\\
   & = I+II,
   \end{align*}
say.  Applying H\"{o}lder's inequality with $\beta > 1$ and using Lemma \ref{lem Vest234}(i) and Lemma \ref{lem stop} we have
    \begin{align*}
    II & \le \left(  \E^x\left[ e^{ \beta  p A_t}\right]\right)^{1/\beta} \PP^x(\tau_k \le t)^{(\beta-1)/\beta}\\
    &  \le \left(  V(x)e^{\Gamma_2 t}\right)^{1/\beta} \left(k^{-1}V(x) e^{\Gamma_0 t}\right)^{(\beta-1)/\beta}\\
    & = V(x)\exp\left\{ \frac{\Gamma_2+(\beta-1) \Gamma_0}{\beta}t\right\} k^{ -(\beta-1)/\beta}.
   \end{align*} 
Since
 \begin{align*}
 \lefteqn{ \sup _{x \in M}\left(\frac{(L^{(k)}+ \beta p \sum_{j=1}^m q_{j,k}X_{j,k} + \beta p Q_k + \frac{\beta^2 p^2}{2}\sum_{j=1}^m q_{j,k}^2)V(x)}{V(x)}\right)} \hspace{25ex} \\
 & =  \sup _{x \in M} \chi_k^2(x) \left(\frac{(L+\beta p Y + \beta p Q + \frac{\beta^2 p^2}{2}R)V(x)}{V(x)}\right)\\
 & \le   \max\left(0,\sup _{x \in M} L_{\beta p} V(x)/V(x)\right)
 \end{align*}
we get a similar estimate for $I$ with the same $\Gamma_2$.  Therefore 
   $$
   \left|T^p_{t,k}g(x) - T^p_t g(x)\right| 
     \le 2 V(x)\exp\left\{ \frac{\Gamma_2+(\beta-1) \Gamma_0}{\beta}t\right\} k^{ -(\beta-1)/\beta}
   $$ 
and the right side converges to 0 as $k \to \infty$ uniformly for $0 \le t \le T$ and $x \in D_{k_0}$.  \qed

\begin{proposition} \label{prop FK phi}  Assume {\bf (0,1,2)}.  For $\phi \in C^\infty_c(M)$ the mapping $(t,x) \to T_t\phi(x)$ satisfies
   \begin{equation} \label{FK phi}
    \left( L +pY + pQ + \frac{p^2}{2}R- \frac{\partial}{\partial t}\right) T_t^p\phi(x) = 0
    \end{equation}
for $(t,x) \in (0,\infty) \times M$ in the sense of distributions.   
\end{proposition}

\n{\bf Proof.}  Write $\widehat{M} = (0,\infty) \times M$ and recall $L_p = L+pY+pQ+\frac{p^2}{2}R$.  It suffices to show that 
   \begin{equation} \label{FK phi2}
   \iint_{\widehat{M}} T_t^p\phi(x)\left(L_p  - \frac{\partial}{\partial t}\right)^\ast \psi(t,x)dtdx = 0
   \end{equation}
for all $\psi \in C^\infty_c(\widehat{M})$, where $dx$ denotes Riemannian measure on $M$ and $(L_p  - \frac{\partial}{\partial t})^\ast$ denotes the adjoint of $L_p  - \frac{\partial}{\partial t}$ with respect to $L^2(\widehat{M}, dtdx)$.  Let $K = \mbox{supp}(\psi)$ and note that 
   $$
   C:= \iint_{\widehat{M}}\left|\left(L_p  - \frac{\partial}{\partial t}\right)^\ast \psi(t,x)\right|dtdx < \infty.
   $$    
    By Lemma \ref{lem FK k} we have 
 $$
   \left(\chi_k^2L_p  - \frac{\partial}{\partial t}\right)T_{t,k}^p\phi(x) = 0     ,
  $$
so that 
    $$
    \iint_{\widehat{M}} T_{t,k}^p\phi(x)\left(\chi_k^2L_p  - \frac{\partial}{\partial t}\right)^\ast \psi(t,x)dtdx = 0.
    $$
There exists $k_0$ such that $K = \mbox{supp}(\psi) \subset (0,\infty)\times D_{k_0}$.    For $k \ge k_0$ we have $\chi_k^2 = 1$ on $D_{k_0}$ and so 
$$
    \iint_{\widehat{M}} T_{t,k}^p\phi(x)\left(L_p  - \frac{\partial}{\partial t}\right)^\ast \psi(t,x)dtdx = 0.
    $$
Then 
   \begin{eqnarray*}
  \lefteqn{ \left|  \iint_{\widehat{M}} T_t^p \phi(x)\left(L_p   - \frac{\partial}{\partial t}\right)^\ast \psi(t,x)dtdx\right|} \hspace{10ex} \\
   & \le & \left|  \iint_{\widehat{M}}\Bigl(T_t^p\phi(x)- T_{t,k}^p\phi(x)\Bigr)\left(L_p  - \frac{\partial}{\partial t}\right)^\ast \psi(t,x)dtdx\right| \\
    & \le & \sup_{(t,x) \in K}\Bigl|T_t^p\phi(x)- T_{t,k}^p\phi(x)\Bigr|  \iint_{\widehat{M}}\left|\left(L_p  - \frac{\partial}{\partial t}\right)^\ast \psi(t,x)\right| dtdx\\
    & = & C \sup_{(t,x) \in K}\left|T_t^p \phi(x)- T_{t,k}^p\phi(x)\right|,
  \end{eqnarray*}
and \eqref{FK phi2} now follows by Lemma \ref{lem kinf}.   \qed

\s

The next task is to obtain a similar result for $u(t,x) = T_t^pg(x)$ with $g \in B_V$.

\begin{lemma} \label{lem FK adj}  Assume {\bf (0,1,2)}.  Suppose $\psi \in C^\infty_c([0,\infty) \times M)$.  Then
    \begin{equation} \label{FK adj}
  \E^x\left[\psi(t,x_t)e^{ pA_t}\right]- \psi(0,x)  =  \E^x \int_0^t \left(L+pY+pQ+ \frac{p^2}{2}R+\frac{\partial }{\partial s}\right) \psi(s,x_s)e^{p A_s}ds 
   \end{equation}
   for all $t \ge 0$ and $x \in M$
\end{lemma}

\n{\bf Proof.} Define 
  $
  M_t = \phi(t,x_t)e^{p A_t}.
  $
It\^{o}'s formula gives 
  \begin{align*}
   M_t-\phi(0,x) &  =  \int_0^t \left(L+pY+p Q+\frac{p^2}{2}R+\frac{\partial }{\partial s}\right) \phi(s,x_s)e^{p A_s}ds \\
   & \quad + \sum_{j=1}^m\int_0^t \big(X_j \phi(s,x_s)+ p q_j(x_s)\phi(s,x_s)\bigr)e^{p A_s}dW_j(s). 
  \end{align*}
The stochastic integrals are local martingales, so we obtain 
  \begin{equation} \label{FK adj n}
  \E^x[ M_{t \wedge \tau_k}]- \phi(0,x)  =  \E^x \int_0^{t \wedge \tau_k}\left(L+pY+pQ+\frac{p^2}{2}R+\frac{\partial }{\partial s}\right) \phi(s,x_s)e^{p A_s}ds 
   \end{equation}
with stopping times $\tau_k \to \infty$.  
Since $\phi$ has compact support in $[0,\infty) \times M$ the integrand on the right of \eqref{FK adj n} is bounded by $C e^{pA_s}$ for some constant $C$ and 
    $$
    \E^x \int_0^t C e^{p A_s}ds = C \int_0^t (T_s^p \mathbf{1})(x) ds
    \le  C \left(\int_0^t e^{\Gamma s}ds\right) V(x) < \infty,
    $$
using Lemma \ref{lem Vest}.   On the left of \eqref{FK adj n} we have $|M_{t \wedge \tau_k}| \le \widetilde{C} e^{p A_{t \wedge \tau_k}}$ for some constant $\widetilde{C}$ and so 
  $$
  \E^x |M_{t \wedge \tau_k}|^\beta \le \widetilde{C}^\beta \E^x\left[e^{\beta  p A_{t \wedge \tau_k}}\right] \le \widetilde{C}^\beta e^{\Gamma_2 t}V(x)
  $$      
using Lemma \ref{lem Vest234}(i).  Thus we can let $k \to \infty$ in \eqref{FK adj n} and obtain \eqref{FK adj}.
\qed

\s

At this point we prepare to use the hypoellipticity assumption {\bf (H)}.   The next lemma extends the assumption from the operator $L = \frac{1}{2} \sum_{j=1}^m X_j^2+X_0$ to the operator $L+pY =   \frac{1}{2} \sum_{j=1}^m X_j^2+(X_0 + pY)$ where $Y(x) = \sum_{j=1}^m q_j(x)X_j(x))$.

\begin{lemma} \label{lem span} Suppose $\widehat{X}_0(x) = X_0(x)+ \sum_{j=1}^m f_j(x)X_j(x)$ for smooth functions $f_1,\ldots,f_m$.  Then 
    \begin{equation} \label{span}
   {\rm span}\{ U(x): U \in  {\cal L}_0(X_0;X_1,\ldots,X_m)\} \subset {\rm span}\{U(x): U \in {\cal L}_0(\widehat{X}_0;X_1,\ldots,X_m)\} 
   \end{equation}
for all $x \in M$.
 \end{lemma}    

\n{\bf Proof.}  
Let $\Gamma(M)$ be the space of all smooth vector fields on $M$. Following Hairer \cite{Hai}, define a collection of vector fields ${\cal X}_k$ by 
 $${\cal X}_0 = \{X_j: j >0\}, \qquad {\cal X}_{k+1} = {\cal X}_k \cup \{[U,X_j]: U \in {\cal X}_k,\, j \ge 0\},
 $$
and then ${\cal X}_k(x) = \mbox{span}\{U(x): U \in {\cal X}_k\}$.  And similarly $\widehat{\cal X}_k$ when $X_0$ is replaced by $\widehat{X}_0$.  Let ${\cal M}(\widehat{\cal X}_k)$ be the $C^\infty(M)$ submodule of $\Gamma(M)$ generated by $\widehat{\cal X}_k.$   An elementary proof using induction on $k$ gives the inclusion 
     $$
     {\cal X}_k \subset {\cal M}(\widehat{\cal X}_k) \qquad \mbox{for all } k \ge 0.
     $$
It follows that 
 $$
    \mbox{span}\{U(x): U \in {\cal X}_k\} \subset \mbox{span}\{U(x): U \in \widehat{\cal X}_k\} \qquad \mbox{for all } k \ge 0
     $$
and then letting $k \to \infty$ gives \eqref{span}.  \qed

\begin{proposition} \label{prop kernel}  Assume {\bf (0,1,2)} and {\bf (H)}.  There exists a $C^\infty$ function $\rho_t(x,y): (0,\infty) \times M \times M \to \R$
 such that 
      \begin{equation} \label{kernel}
      T_t^p g(x) = \int_M \rho_t(x,y)g(y)dy
     \end{equation}
for all $g \in B_V$. Further 
     $$
     \left(L_{p,x} - \frac{\partial}{\partial t}\right)\rho_t(x,y) = 0 \quad \mbox{ and } \quad  \left(L_{p,y}^\ast - \frac{\partial}{\partial t}\right)\rho_t(x,y) = 0,
     $$
where $L_{p,x}$ denotes the action of $L_p=L+pY+pQ+\frac{p^2}{2}R$ on functions of $x$, and similarly for $L_{p,y}$.  
 \end{proposition}

\n{\bf Proof.}   We extend Theorem 3 of Ichihara and Kunita \cite{IK} from the case involving the semigroup of the generator $L$ to the more general case of the semigroup $\{T_t^p: t \ge 0\}$, using \eqref{FK adj} in place of \cite[eqn (4.3)]{IK}.  Write $T_t^p$ in integral operator form:  
        $$
        T_t^p g(x) = \int_M R_t(x,dy)g(y).
        $$       
Putting $\psi \in C^\infty_c((0,\infty) \times M)$ in \eqref{FK adj} and letting $t \to \infty$ gives
 $$
     \int_0^\infty \int_M R_s(x,dy)\left(L_p+\frac{\partial }{\partial s}\right) \psi(s,y) ds = 0
   $$
or equivalently
      $$
      \left(L_p+\frac{\partial }{\partial s}\right)^\ast_{(s,y)}R_s(x,dy) = 0
      $$
in the sense of distributions.  By assumption {\bf (H)} and Lemma \ref{lem span} the operator $L_p +\frac{\partial}{\partial s}$ is hypoelliptic, and so $\left(L_p +\frac{\partial}{\partial s}\right)^\ast$ is also hypoelliptic.  It follows that for each $x$ we have $R_s(x,dy) = \rho_s(x,y)dy$ for some $\rho_s(x,y)$ which is smooth in $(s,y)$.  The rest of the argument in the proof of \cite[Theorem 3]{IK} goes through with $L$ replaced by $L_p $, and with the lemma of Blagore\u{s}\u{c}ensky and Freidlin \cite[Lemma 4.1]{IK} replaced by Proposition \ref{prop FK phi} above.      \qed      
 
\s

 Without some integrability conditions on $\rho_t(x,y)$ and its partial derivatives we cannot immediately deduce from \eqref{kernel} that $T_t^pg$ is continuous when $g \in B_V$, or that 
     $$
    \left(L_p-\frac{\partial}{\partial t}\right) T_t^pg(x)
     =  \int_M \left( \bigl(L_p-\frac{\partial}{\partial t}\bigr)\rho_t(x,y)\right)g(y)dy = 0.
     $$ 
The following result uses just the continuity of $\rho_t(x,y)$ together with Lemma \ref{lem Vest12}.

\begin{lemma} \label{lem comp} Denote $C_V = B_V \cap C(M)$.  Assume {\bf (0,1,2)} and the parabolic H\"{o}rmander condition {\bf (H)}.  For each $t> 0$ the operator $T_t^p$ is a compact operator on $B_V$ and maps into $C_V$.
\end{lemma}  
 
\n{\bf Proof.}  We follow the proof of compactness in Rey-Bellet \cite[Theorem 8.9]{RB06}.  Recall the compact sets $D_k = \{x : V(x) \le k\}$.  Write
  \begin{align*}
     T^p_t & = (\mathbf{1}_{D_k} + \mathbf{1}_{D_k^c})T^p_{t/2}(\mathbf{1}_{D_k} + \mathbf{1}_{D_k^c})T^p_{t/2} \\
        & = \mathbf{1}_{D_k}T^p_{t/2}\mathbf{1}_{D_k}T^p_{t/2} + \mathbf{1}_{D_k}T^p_{t/2}\mathbf{1}_{D_k^c}T^p_{t/2} + \mathbf{1}_{D_k^c}T^p_{t/2}\mathbf{1}_{D_k}T^p_{t/2} +\mathbf{1}_{D_k^c}T^p_{t/2} \mathbf{1}_{D_k^c}T^p_{t/2},
   \end{align*}    
where $\mathbf{1}_{D_k}$ denotes multiplication by the indicator function of $D_k$, and similarly for $\mathbf{1}_{D_k^c}$.  By Lemma \ref{lem Vest12} we have $T^p_{t/2}V(x)/V(x) \to 0$ as $x \to \infty$, and so the operator norm $\|\mathbf{1}_{D_k^c}T^p_{t/2}\|_{B_V\to B_V} \to 0$ as $k \to \infty$.  Therefore 
   $$
   \|T^p_t - \mathbf{1}_{D_k}T^p_{t/2}\mathbf{1}_{D_k}T^p_{t/2}\|_{B_V\to B_V} \le 3 \|T^p_{t/2}\|_V \|\mathbf{1}_{D_k^c} T^p_{t/2}\|_{B_V\to B_V} \to 0
   $$
as $k \to \infty$.  The uniform continuity of $(x,y) \to \rho_{t/2}(x,y)$ on the compact set $D_k \times D_k$ together with the theorem of Arzel\'{a} and Ascoli implies that $\mathbf{1}_{D_k}T^p_{t/2}\mathbf{1}_{D_k}$ is a compact operator.   Therefore $T^p_t$ is the limit of compact operators $\bigl(\mathbf{1}_{D_k}T^p_{t/2}\mathbf{1}_{D_k}\bigr)T^p_{t/2}$, and so $T^p_t$ is compact. 

Now suppose $g \in B_V$ and write $h_k(x) = \mathbf{1}_{D_k}T^p_{t/2} \mathbf{1}_{D_k}T^p_{t/2}g(x)$.  Then $h_k$ is continuous on $D_k$ and $h_k$ converges to $T^p_tg$ uniformly on compact sets as $k \to \infty$.  Therefore $T^p_tg$ is continuous.  \qed

\begin{lemma} \label{lem gphin}  Assume {\bf (0,1,3)}.   Given $g \in C_V$ there exists a sequence $\phi_k \in C^\infty_c(M)$ such that $T^p_t\phi_k(x) \to T^p_tg(x)$ as $k \to \infty$ uniformly on compact subsets of $[0,\infty) \times M$.
\end{lemma}  

\n{\bf Proof.}  Without loss of generality we may assume $\|g\|_V \le 1$. Recall the sets $D_k = \{x\in M: V(x)   \le k\}$ and the smooth functions $\chi_k$.  Let $g_k(x) = \chi_k(x)g(x)$.  Since $g_k$ is continuous with $\mbox{supp}(g_k) \subset \mbox{int}(D_{k+1})$ there exists $\phi_k \in C^\infty_c(M)$ such that $|g_k(x)-\phi_k(x)| \le 1/k$ for all $x \in M$.  For $x \in D_k$ we have $|g(x)-\phi_k(x)| \le 1/k$ and for $x \not\in D_k$ we have $|g(x) - \phi_k(x)| \le |g_k(x)-\phi_k(x)|+|g(x)-g_k(x)| \le 1/k+ V(x)$.  Together
    $$
    |g(x)-\phi_k(x)| \le \frac{1}{k} + V(x)1_{V(x) > k}   \le \frac{1}{k}+  \frac{1}{k^{\beta-1}}V^\beta(x)
    $$
for any $\beta > 1$.  Choose $\beta > 1$ so that growth condition {\bf (3)} is satisfied.  By Lemmas \ref{lem Vest} and \ref{lem Vest234}(ii) we have 
   $$
   |T_t^pg(x) - T_t^p\phi_k(x)| \le \frac{1}{k}T_t^p{\bf 1}(x)+ \frac{1}{k^{\beta-1}}T_t^pV^\beta(x) \le \frac{1}{k}e^{\Gamma t}V(x)+ \frac{1}{k^{\beta-1}} e^{\Gamma_3 t}V^\beta(x)
   $$
for all $t \ge 0$, and the result follows directly.  \qed

\begin{proposition} \label{prop FK g} Assume {\bf (0,1,2,3)} and {\bf (H)}.  For $g \in B_V$ the mapping $(t,x) \to T_t^p g(x)$ is $C^\infty$ on $(0,\infty )\times M$ and satisfies
   \begin{equation} \label{FK g}
    \left( L + pY+pQ+ \frac{p^2}{2}R - \frac{\partial}{\partial t}\right) T_t^p g(x) = 0, \quad \mbox{for all }t > 0, x \in M. 
    \end{equation}
   \end{proposition}

\n{\bf Proof of Proposition \ref{prop FK g}.}  By Lemma \ref{lem comp} for any $t_0 > 0$ the function $T^p_{t_0}g \in C_V$, and then $T_t^pg = T^p_{t-t_0}T^p_{t_0}g$ for $t > t_0$.  So without loss of generality we may assume $g \in C_V$.  Let $\phi_k \in C^\infty_c(M)$ be the sequence in Lemma \ref{lem gphin}.   We adapt the proof of Proposition \ref{prop FK phi}. 
Since $L_p-\frac{\partial}{\partial t}$ is hypoelliptic on $\widehat{M} = (0,\infty)\times M$ it suffices to show that 
   \begin{equation} \label{FK g2}
   \iint_{\widehat{M}} T_t^p g\left(L_p - \frac{\partial}{\partial t}\right)^\ast \psi(t,x)dtdx = 0
   \end{equation}
for all $\psi(t,x) \in C^\infty_c(\widehat{M})$.  Let $K = \mbox{supp}(\psi)$ and note that 
   $$C:= \iint_{\widehat{M}}\left|\left(L_p - \frac{\partial}{\partial t}\right)^\ast \psi(t,x)\right|dtdx < \infty.
   $$     By Proposition \ref{prop FK phi} we have 
 $$
   \left(L_p - \frac{\partial}{\partial t}\right) T^p_t \phi_k(x) = 0
  $$
for each $k$, and then 
   \begin{eqnarray*}
  \lefteqn{ \left|  \iint_{\widehat{M}} T_t^p g(x)\left(L_p  - \frac{\partial}{\partial t}\right)^\ast \psi(t,x)dtdx\right|} \hspace{10ex} \\
   & \le & \left|  \iint_{\widehat{M}}\Bigl(T_t^p g(x)- T_t^p \phi_n(x)\Bigr)\left(L_p - \frac{\partial}{\partial t}\right)^\ast \psi(t,x)dtdx\right| \\
    & \le & \sup_{(t,x) \in K}\Bigl|T_t^p g(x)- T_t^p \phi_k(x)\Bigr|  \iint_{\widehat{M}}\left|\left(L_p - \frac{\partial}{\partial t}\right)^\ast \psi(t,x)\right| dtdx\\
    & = & C \sup_{(t,x) \in K}\left|T_t^p g(x)- T_t^p\phi_k(x)\right|,
  \end{eqnarray*}
and \eqref{FK g} now follows by Lemma \ref{lem gphin}.   \qed

\subsection{Positivity}

The proof of Theorem \ref{thm KR} is based on the Kre\u{\i}n-Rutman theorem for compact positive operators on a Banach space.   
Positivity in $B_V$ is defined in terms of the total cone ${\cal K} = \{g \in B_V: g(x) \ge 0 \mbox{ for all }x \in M\}$.  It is immediate from the definition of $T_t^p$ that $T_t^pg \in {\cal K}$ for all $g \in {\cal K}$.  Thus $T_t^p$ is a positive operator on $B_V$.  
 However we will need a stronger form of positivity, and this is where assumption {\bf (P)} is used.  

\begin{lemma} \label{lem control} Assume {\bf (H,P,G)}.  Suppose $x \in M$ and $t > 0$, and let $U$ be a non-empty open subset of $M$. Then $T_t^p \mathbf{1}_U(x) = \E^x\left[e^{pA_t}1_{x_t \in U}\right] > 0$.
\end{lemma}

\n{\bf Proof.}   This follows immediately from {\bf (P)} since $\PP^x(e^{pA_t} > 0) = 1$.   \qed 
\s

\begin{lemma} \label{lem spr} Assume {\bf (H,P,G)}.   For $t > 0$ the operator $T_t^p$ has strictly positive spectral radius.
\end{lemma}  

\n{\bf Proof.}  Choose $g = \mathbf{1}_U$ where $U \subset M$ is open and non-empty with compact closure.  Lemma \ref{lem control} implies that $T_t^pg(x)> 0$ for all $x \in M$.  Let $c = \inf \{T_t^p g(x): x \in U\} > 0$, so that $T_t^pg(x) \ge cg(x)$ for all $x \in M$.  Iterating gives $(T^p_t)^n g(x) \ge c^n g(x)$ for all $x \in M$, so the spectral radius of $T_t^p$ is at least $c >0$.  \qed

\begin{lemma} \label{lem petite}  Assume {\bf (H,P,G)}.  For all $t \ge 0$ and compact $C \subset M$ there exists a relatively compact non-empty open set $U_1 \subset M$ and $\delta > 0$ such that $\rho_t(x,y) \ge \delta$ for all $x \in C$ and $y \in U_1$.
\end{lemma}

\n{\bf Proof.}   The proof is an easy adaptation of the argument in \cite[Page 223]{MSH02}.  There exist $x_0$ and $x_1$ such that $\rho_{t/2}(x_0,x_1) > 0$.  By continuity of $\rho_{t/2}$ there exists (non-empty) open sets $U$ and $U_1$ and $\delta_1 > 0$ such that $\rho_{t/2}(z,y) \ge \delta_1$ for $z \in U_0$ and $y \in U_1$.  Also $x \mapsto T_{t/2}^p \mathbf{1}_{U_0}(x) = \int_{U_0} \rho_{t/2}(x,z)dz$ is continuous and strictly positive so there exists $\delta_2$ such that $\int_{U_0} \rho_{t/2}(x,z)dz \ge \delta_2$ for all $x \in C$.   Then for $x \in C$ and $y \in U_1$ we have 
   $$
   \rho_t(x,y) \ge \int_{U_0} \rho_{t/2}(x,z)\rho_{t/2}(z,y)dz \ge \delta_1 \int_{U_0} \rho_{t/2}(x,z)dz \ge \delta_1 \delta_2,
   $$
and the proof is complete.  \qed   
   
\s

Lemma \ref{lem control} implies that $T_t^pg(x) > 0$ for all $x$ for all $g \in {\cal K} \setminus \{0\}$.  
However Lemma \ref{lem Vest12} implies that $T^p_tg(x)/V(x) \to 0$ as $x \to \infty$ for all $g \in B_V$.  Therefore for $g \in {\cal K}$ we have $T_t^p g \not\in \overset{\circ}{\cal K}$, and so $T_t$ is not strongly positive for the cone ${\cal K} \subset B_V$.  In the proof of Theorem \ref{thm KR} we use the version of the Kre\u{\i}n-Rutman theorem given in Deimling \cite[Thm 19.2, page 226]{Deim}, but not the version requiring strong positivity given in \cite[Thm 19.3, page 228]{Deim}.

\subsection{Proof of Theorem \ref{thm KR}}

The proof follows closely the proof of \cite[Lemma 6.6]{FS20}.   Let $\kappa_p$ denote the spectral radius of the time $1$ operator $T_1^p$.  Lemma \ref{lem spr} shows $\kappa_p > 0$.  The existence of the function $\phi_p \in {\cal K} \cap B_V$ satisfying $T_1^p\phi_p(x) = \kappa_p \phi_p(x)$ is due to the Kre\u{\i}n-Rutman theorem, see for example Deimling \cite[Thm 19.2, p 226]{Deim}.
  Since $\phi(x) = \kappa_p^{-1}T_1^p\phi_p(x)$ then $\phi_p \in C^\infty(M)$ by Theorem \ref{thm Tt}(ii).  
   Since $\phi_p \in {\cal K}$ is not identically zero there exists $x_0$ such that $\phi_p(x_0) > 0$.  Since $\phi_p$ is continuous there exists a non-empty open set $U \subset M$ such that $\phi_p \ge \delta \mathbf{1}_U$ for some $\delta > 0$, and then $\phi_p(x) = \kappa_p^{-1}T_1^p\phi_p(x) \ge \delta \kappa_p^{-1}T_1^p \mathbf{1}_U(x) > 0$ for all $x \in M$.  
   
   Define 
      $$
      \widehat{T}g(x) = \frac{1}{\kappa_p \phi_p(x)}T_1^p(\phi_p g)(x)
     $$
so that $\widehat{T}$ is a Markov operator with 
transition density $\widehat{\rho}(x,y) = \dfrac{1}{\kappa_p \phi_p(x)}\rho_1(x,y) \phi_p(y)$.  (This corresponds to the {\it twisted kernel} in \cite[Eqn (11)]{KM03}.)  Define $\widehat{V}(x) = V(x)/\phi_p(x) > 0$.   By Lemma \ref{lem Vest12} 
     $$
     \frac{\phi_p(x)}{V(x)} = \frac{T_1^p\phi_p(x)}{\kappa_p V(x)} \le \frac{\|\phi_p\|_V}{\kappa_p} \frac{T_1^pV(x)}{V(x)} \to 0 \quad \mbox{ as } x \to \infty
     $$
so that $\widehat{V}(x) \to \infty$ as $x \to \infty$ and in particular $\widehat{V}$ is bounded away from 0.  Also 
    \begin{equation} \label{V4}
      \frac{\widehat{T}\widehat{V}(x)}{\widehat{V}(x)} = \frac{T_1^pV(x)}{\kappa_p V(x)} \to 0 \quad \mbox{ as } x \to \infty.
     \end{equation}    
Lemma \ref{lem petite} together with the strict positivity and continuity of $\phi_p$ implies that every compact set $C \subset M$ is petite for the Markov operator $\widehat{T}$.  The estimate \eqref{V4} implies that $\widehat{V}$ satisfies condition (V4) of Meyn and Tweedie \cite{MTbook}.  It follows by \cite[Theorem 16.1.2]{MTbook} that the discrete time Markov chain with operator $\widehat{T}$ is $\widehat{V}$-exponentially ergodic.  

There exists a unique $\widehat{T}$-invariant probability measure $\pi$, say, and $\int_M \widehat{V} d\pi < \infty$.  For non-empty open $U \subset M$ we have $T_1^p\mathbf{1}_U(x) > 0$ for all $x$, so that $\widehat{T} \mathbf{1}_U(x) > 0$ for all $x$, and hence $\pi(U)  > 0 $.  The smoothness of $\rho_1(x,y)$, and hence of $\hat{\rho}(x,y)$, imply that $\pi$ has a smooth density with respect to Lebesgue measure on $M$.  Further there exist constants $C_1$ and $\eta < 1$ 
such that
    \begin{equation} \label{gehat}
  \left| \widehat{T}^n \widehat{g}(x) -  \int_M \widehat{g} d\pi \right|  \le C_1 \eta^n \widehat{V}(x)
  \end{equation}
for all $x \in M$ and all $n \ge 0$ whenever $|\widehat{g}(x)| \le \widehat{V}(x)$ for all $x \in M$.   Define the measure $\widehat{\nu}_p$ by $\dfrac{d\widehat{\nu}_p}{d\pi}(x) = \dfrac{1}{\phi_p(x)}$.  Then
    $\int_M V d\widehat{\nu}_p =  \int_M \widehat{V} d\pi < \infty$ and $\widehat{\nu}_p(U) > 0$ for all non-empty open $U \subset M$.

Suppose $g \in B_V$.  Putting $\widehat{g}(x) = g(x)/\phi_p(x)$ in \eqref{gehat} and multiplying by $\phi_p(x)$ gives
    \begin{equation} \label{gen}
  \left| \frac{1}{(\kappa_p)^n}T_n^p g(x) -  \phi_p(x) \int_M g d\widehat{\nu} \right|  \le C_1 \|g\|_V \eta^n V(x)
  \end{equation}
for all $x \in M$ and $n \ge 0$. 

Suppose now $g = g_1+ig_2$ in the complexification $\mathbb{C}B_V$ of $B_V$ is a (nontrivial) eigenfunction of $T_1^p$ with eigenvalue $\tilde{\kappa}\in \mathbb{C}$, say.  Applying \eqref{gen} to $g_1$ and $g_2$ and combining the results gives
      $$
      \left(\frac{\tilde{\kappa}}{\kappa_p}\right)^n g(x) \to  \phi_p(x) \int_M g\,d\widehat{\nu}
      $$
as $n \to \infty$.   
If $\int g\,d\widehat{\nu} = 0$ then $|\tilde{\kappa}| < \kappa$, and if $\int g\,d\widehat{\nu} \neq 0$ then $\tilde{\kappa} = \kappa$ and $g$ is a (complex) multiple of $\phi_p$.  It follows that the positive eigenfunction $\phi_p$ is unique up to a multiplication by a positive scalar.   

Applying the argument above to the operator $T_{1/m}^p$ in place of $T_1^p$ gives a strictly positive eigenfunction $\phi_{p,m} \in B_V$ and eigenvalue $\kappa_{p,m}$. Since $T_1^p\phi_{p,m} = (T_{1/m}^p)^m \phi_{p,m} = (\kappa_{p,m})^m \phi_{p,m}$ we deduce that $\phi_{p,m} = \phi_p$ (up to a scalar multiple) and $(\kappa_{p,m})^m = \kappa_p$.  So we can take $\phi_{p,m} = \phi_p$ and $\kappa_{p,m} = (\kappa_p)^{1/m}$.  For any rational time $t = \ell/m $ we have $T_{\ell/m}^p\phi_p = (T_{1/m}^p)^{\ell}\phi_p = (\kappa_p)^{\ell/m}\phi_p$.  Further, for any time $t >0$, let $t_k$ be a sequence of rational times with $t_k \to t$.  Then by Theorem \ref{thm Tt}(ii),
      $$
      T_t^p \phi_p(x) =  \lim_{k \to \infty} T_{t_k}^p\phi_p(x) = \lim_{k \to \infty} (\kappa_p)^{t_k}\phi_p(x) = \kappa^t \phi_p(x),
      $$ 
so that $\phi_p$ is an eigenfunction of $T_t$ with eigenvalue $\kappa^t$.

Next we extend \eqref{gen} from discrete time $n$ to continuous time $t$.   
 Note that \eqref{gen} is equivalent to the statement
    $$
    \left\|\frac{1}{(\kappa_p)^n} T^n g - \phi_p \int g d \widehat{\nu}_p\right\|_V \le C_1 \|g\|_V \eta^n.
    $$
Writing $t = n+s$ for integer $n \ge 0$ and $0 \le s \le 1$ we have
 \begin{align*}
    \left\|\frac{1}{(\kappa_p)^t} T_t^p g - \phi_p \int_M g d\widehat{\nu}_p \right\|_V 
    & =  \left\|\frac{1}{(\kappa_p)^s}T_s^p\left(\frac{1}{(\kappa_p)^n} T_n^p g - \phi_p \int_M g d\widehat{\nu}_p \right) \right\|_V \\
     & \le   \frac{\|T_s^p\|_V}{(\kappa_p)^s}\left\|\left(\frac{1}{(\kappa_p)^n} T_n^p g - \phi_p \int g d\widehat{\nu}_p \right) \right\|_V \\   
            &  \le  \frac{e^{\Gamma s}}{(\kappa_p)^s} C_1 \|g\|_V \eta^n =  \left(\frac{e^{\Gamma }}{\eta \kappa_p}\right)^s C_1\|g\|_V \eta^t 
      \end{align*}   
and this gives \eqref{get} with $C = C_1\max(e^{\Gamma}/(\eta \kappa_p),1)$.

   Finally the proof for any eigenfunction $g$ of the complexification of $T_t^p$ on $\mathbb{C}B_V$ is a repeat of the argument above, using \eqref{get} in place of \eqref{gen}.   \qed

\subsection{Proof of Corollary \ref{cor Lam}}

 \n{\bf Proof of Corollary \ref{cor Lam}(i).}   By Theorem \ref{thm KR} we have 
     $$
     \frac{1}{\kappa_p^t} T_t^p\mathbf{1}(x) \to \phi_p(x) \int_M \mathbf{1}\,d\widehat{\nu}_p >0
     $$ 
so that
    $$   
   \log  T_t^p \mathbf{1}(x) - t \log \kappa_p \to \log \left( \phi_p(x)\int_M \mathbf{1}\,d\widehat{\mu}_p \right) \in (-\infty,\infty)
    $$
and thus 
    $$
   \frac{1}{t} \log \E^x \left[e^{pA_t}\right] =  \frac{1}{t} \log T_t^p \mathbf{1}(x) \to \log \kappa_p
    $$
as $t \to \infty$ for all $x \in M$. \qed

\s

 \n{\bf Proof of Corollary \ref{cor Lam}(ii).}   For $t > 0$ we have $T_t ^p\phi_p(x) = e^{\Lambda(p)t}\phi_p(x)$, and the result follows directly from Theorem \ref{thm Tt}(ii).  \qed
  
\s

 \n{\bf Proof of Corollary \ref{cor Lam}(iii).}   Let $M_t = e^{-\widetilde{\Lambda}t} g(x_t) e^{p A_t}$. 
An application of It\^{o}'s formula shows that $M_t$ is a local martingale, so there are stopping times $\tau_n \to \infty$ such that $\E^x[M_{t \wedge \tau_n}] = g(x)$.   The bound
 \begin{align*}
  \E^x\left[\left(V(x_{t \wedge \tau_n})e^{pA_{t \wedge \tau_n}}\right)^{(1+\beta)/2}\right] 
   & \le \frac{1}{2} \left(\E^x\left[V(x_{t \wedge \tau_n})e^{\beta  pA_{t \wedge \tau_n}}\right] + \E^x\left[V^\beta(x_{t \wedge \tau_n})e^{ pA_{t \wedge \tau_n}}\right]\right) \\
   & \le   \frac{1}{2}\left(e^{\Gamma_2 t}V(x) + e^{\Gamma_3 t}V^\beta (x)\right) 
  \end{align*}
for some $\beta >1$ from Lemma \ref{lem Vest234} implies that the random variables $V(x_{t \wedge \tau_n}) e^{ pA_{t \wedge \tau_n}} $ are uniformly $\PP^x$-integrable.  Since $g \in B_V$ it follows that the $M_{t \wedge \tau_n}$ are uniformly $\PP^x$-integrable and so 
    $$g(x) = \E^x[M_{t \wedge \tau_n}] \to \E^x[M_t] = e^{-\widetilde{\Lambda}t} T_t g(x).
    $$
Therefore $g$ is an eigenfunction of $T_t^p$ with eigenvalue $e^{t \widetilde{\Lambda}}$, and the result follows from the uniqueness of non-negative eigenfunctions in Theorem \ref{thm KR}.  \qed

\begin{remark}  \label{rem pos} Suppose that instead of the assumption in Corollary \ref{cor Lam}(iii) we have $g \in B_V \cap C^2(M)$ is non-negative, not identically 0, and satisfies $L_p g \ge \widetilde{\Lambda} g$.  Then $M_t$ is a local submartingale and so there are stopping times $\tau_n \to \infty$ such that $\E^x[M_{t \wedge \tau_n}] \ge g(x)$.  As above, the $M_{t \wedge \tau_n}$ are uniformly $\PP^x$-integrable and so    $$g(x) \le \E^x[M_{t \wedge \tau_n}] \to \E^x[M_t] = e^{-\widetilde{\Lambda}t} T_t^pg(x).
    $$
That is, $T_t^pg \ge e^{\widetilde{\Lambda}t} g$.  Therefore the spectral radius of $T_1^p$ is at least $e^{\widetilde{\Lambda}}$ and so $\Lambda(p) = \log \kappa_p \ge \widetilde{\Lambda}$. 
\end{remark}

\section{Proofs for Section \ref{sec vary}} \label{sec vary proof}

We will use the analytic perturbation theory of Kato \cite{Kato} for operators acting on a complex vector space.  The complexification of the Banach lattice $B_V$ is $\mathbb{C}B_V = \{g_1+ig_2: g_1,g_2 \in B_V\}$
with the norm 
     $\|g_1+ig_2\|_{\mathbb{C}V} = \big\| |g_1+ig_2| \big\|_V $, 
     see for example \cite[page 330]{DJT}.  
Let $\|T\|_{\mathbb{C}V}$ denote the operator norm of a bounded linear operator $T:\mathbb{C}B_V \to \mathbb{C}B_V$.  If $T:B_V \to B_V$ then the extension of $T$ to $\mathbb{C}B_V$ satisfies $\|T\|_{\mathbb{C}V} \le \sqrt{2} \|T\|_V$.        
 
We may extend the definition of the semigroup of operators $\{T_t^p: t \ge 0\}$ acting on $B_V$ for real $p$ to the semigroup of operators $\{T_t^z: t \ge 0\}$ acting on $\mathbb{C}B_V$ for complex $z$ using the formula 
    \begin{equation} \label{Tz}
      T_t^z(g_1+ig_2)(x) = \E^x\left[\bigl(g_1(x_t)+ig_2(x_t)\bigr)e^{zA_t}\right], 
      \end{equation}
noting that when $z$ is real then \eqref{Tz} is simply the complexification of the operator defined by \eqref{Ttp}.

\begin{lemma} \label{lem T anal}  Assume {\bf (H,P,G)} valid for two values $p = p_0 \pm \delta$ and the same $V$ for some $\delta > 0$.   Fix $t_0 > 0$.  For $n \ge 0$ define
    $$
    T_n g(x) = \frac{1}{n!} \E^x\left[g(x_{t_0})e^{p_0A_{t_0}} (A_{t_0})^n\right].
    $$

(i)  For each $n \ge 0$ the map $T_n$ is a bounded linear operator of $B_V$ into itself.

(ii)  If $p_0-\delta < p < p_0+\delta$ then  
    \begin{equation} \label{T anal}
 \sum_{n=0}^\infty|(p-p_0)^n| \|T_n\|_V < \infty \quad \mbox{ and } \quad T_{t_0}^p  =  \sum_{n=0}^\infty  (p-p_0)^n T_n.
  \end{equation}

(iii)  The family $\{T_{t_0}^p: p_0-\delta < p  < p_0+\delta\}$ of operators on $B_V$ has a natural extension to an analytic family $\{T_{t_0}^z: z \in D(p_0,\delta)\}$ of operators on the complexification $\mathbb{C}B_V$ of $B_V$.  Here $D(p_0,\delta) = \{z \in \mathbb{C}: |z - p_0| < \delta\}$ denotes the open ball with center $p_0$ and radius $\delta$ in $\mathbb{C}$. \end{lemma}

\n{\bf Proof.} The convexity of \eqref{LpV} with respect to $p$ implies that {\bf (G)} is satisfied for all $p \in [p_0-\delta,p_0+\delta]$.  Further
   $$
   \sup_{|p-p_0|\le \delta} \sup_{x \in \R} \frac{L_pV(x)}{V(x)} = \max\left(\sup_{x \in \R} \frac{L_{p_0-\delta}V(x)}{V(x)},\sup_{x \in \R} \frac{L_{p_0+\delta}V(x)}{V(x)} \right) < \infty,
   $$ 
and so there exists $\Gamma$ such that $L_pV(x)/V(x) \le \Gamma$ for all $x \in M$ and $p_0-\delta \le p \le p_0+\delta$.  By Lemma \ref{lem Vest} we have 
  \begin{equation} \label{claim}
      T_t^pV(x) \le e^{\Gamma t}V(x)
    \end{equation}
for all $x \in M$ and $p_0-\delta \le p \le p_0+\delta$.  For $n \ge 0$ and $g \in B_V$,
  \begin{align}
  \frac{1}{n!}\E^x\left|g(x_{t_0})e^{p_0A_{t_0}}(A_{t_0})^n\right|
    & \le \frac{1}{n!}\|g\|_V\E^x\left[V(x_{t_0})e^{p_0A_{t_0}} |A_{t_0}|^n\right] \nonumber \\
     & \le \frac{1}{\delta^n}\|g\|_V\E^x\left[V(x_{t_0})e^{p_0A_{t_0}}e^{\delta|A_{t_0}|}\right] \nonumber \\
      & \le \frac{1}{\delta^n}\|g\|_V\Bigl(T_{t_0}^{p_0+\delta} V(x)+ T_{t_0}^{p_0-\delta}V(x) \Bigr)  \nonumber \\
       & \le \frac{2e^{\Gamma t_0}}{\delta^n}\|g\|_V V(x). \label{T anal2}
     \end{align}
It follows that $T_n$ is a well-defined operator on $B_V$ and $\|T_n\|_V \le 2 e^{\Gamma t_0}/\delta^n$.  The first part of \eqref{T anal} follows immediately, and the second part of \eqref{T anal} follows by applying the dominated convergence theorem to the expression
     $$
     g(x_{t_0})e^{p A_{t_0}} = \sum_{n=0}^\infty \frac{1}{n!} g(x_{t_0})e^{p_0A_{t_0}}(p-p_0)^n (A_{t_0})^n.
     $$
Finally, the analyticity of the operators $\{T_{t_0}^z: |z -p_0| < \delta \}$ is obtained by replacing $p$ with $z$ and $g$ with $g_1+ig_2$ in the proof of the second part of \eqref{T anal}.\qed

\s

For the next stage of the argument fix $t_0 =1$.  
   By Theorem \ref{thm KR} the spectrum $\Sigma(T_1^{p_0})$ of the mapping $T_1^{p_0}$ 
acting on $\mathbb{C}B_V$ has an isolated point $\kappa_{p_0} > 0$, say, with (algebraic) multiplicity 1, and the rest of spectrum of $T_1^{p_0}$ lies within the closed ball $\overline{D}(\eta) \subset \mathbb{C}$ of radius $\eta$ for some $\eta < \kappa_{p_0}$, that is, $\Sigma(T_1^{p_0}) \setminus \{\kappa_{p_0}\} \subset \overline{D}(\eta)$. 

 The next result is an application of Kato's analytic perturbation theory \cite[Theorems VII-1.7 and VII-1.8 (pp 368--370)]{Kato}.  
 
\begin{theorem} \label{thm kato}  Assume {\bf (H,C,G)} valid for two values $p = p_0 \pm \delta$ and the same $V$ for some $\delta > 0$.  There exists $0 < \delta_1 \le \delta$ and analytic functions $z \to \widehat{\kappa}_z \in \mathbb{C}$ and $z \to P^z \in L(\mathbb{C}B_V)$ for $|z-p_0|  <\delta_1$ with $\widehat{\kappa}_{p_0} = \kappa_{p_0}$ so that $\widehat{\kappa}_z$ is an eigenvalue of $T^z_1$ of (algebraic) multiplicity 1 and $P^z$ is the corresponding eigenprojection.  In particular $T_1^z P^z g = \widehat{\kappa}_z P^zg$ for all $g \in \mathbb{C}B_V$.  Moreover, $\delta_1$ may be chosen sufficiently small that $\widehat{\kappa}_z$ is the spectral radius of $T_1^z$ for $|z-p_0| <\delta_1$.
\end{theorem}

Notice that for real $p \in (p_0-\delta_1, p_0+\delta_1)$ we may identify the spectral radius $\widehat{\kappa}_p$ above with the spectral radius $\kappa_p$ in Theorem \ref{thm KR}.  Thus the function $\{\widehat{\kappa}_z : |z-p_0| < \delta_1\}$ is an extension of $\{\kappa_p: |p-p_0| < \delta_1\}$.   For ease of notation we henceforth write $\widehat{\kappa}_z = \kappa_z$ for $|z-p_0|< \delta_1$.



\subsection{Proof of Theorem \ref{thm kato1}}

Since the operator $L_p$ is convex in $p$, then the conditions {\bf(H,P,G)} are satisfied for all $p \in [p_1, p_2]$.  
H\"{o}lder's inequality implies that the mapping $p \to \frac{1}{t}\log \E^x[e^{pA_t}]$ is convex for each $t > 0$, and the convexity is preserved for the pointwise limits as $t \to \infty$.  

For each $p_0 \in (p_1,p_2)$ there exists $\delta >0$ such that $p_1 \le p_0-\delta < p_0+\delta \le p_2$.  So we may apply Theorem \ref{thm kato} and obtain $\delta_1 \le \delta$ and the functions $\kappa_z$ and $P^z$.  
 Restricting to real $p$, we deduce that 
   $p \to \Lambda(p) = \log \kappa_p$ and $p \to \phi_p  \equiv P^p(\mathbf{1}) \in \mathbb{C}B_V$ are real analytic on $(p_0-\delta_1,p_0+\delta_1)$.  
 Since $\Lambda(p)$ is defined uniquely for each $p$ it follows immediately that $\Lambda(p)$ is real analytic on $(p_1,p_2)$.  

 Piecing together real analytic functions $p \to \phi_p$ on intervals $(p_0-\delta_1,p_0+\delta_1)$ to obtain a real analytic function $p \to \phi_p$ on the entire interval $(p_1,p_2)$ requires more work, because each eigenfunction is unique only up to a non-zero multiplicative constant.  
 Thus the family $\{\phi_p : p \in I\}$ on some open interval $I$ could be replaced by the family $\{h(p) \phi_p: p \in I\}$ for any non-zero function $h$. 
  We can impose uniqueness by choosing a point $x_0 \in M$ and then insisting $\phi_p(x_0) = 1$ for all $p$.  That is, we can replace a real analytic family $\{\phi_p: p \in I\}$ with the real analytic family $\{\frac{1}{\phi_p(x_0)}\phi_p: p \in I\}$.   
  With the extra requirement on the eigenfunction $\phi_p$ that $\phi_p(x_0) =1$ then $\phi_p$ is uniquely defined for all $p \in (p_1,p_2)$ and the real analytic behavior on each subinterval $(p_0-\delta_2,p_0+\delta_2\}$ implies real analytic behavior on the full interval $(p_1,p_2)$. \qed

\subsection{Proof of Theorem \ref{thm kato2}}

The basic idea is to control the size of the gap in the spectrum of $T_1^z$ for $z$ near to $p_0$.  More specifically we wish to control the size of $(T_1^z - T_1^zP^z)^n$ for large $n$, uniformly for $z$ near $p_0$. 
   Recall $\eta < \kappa_{p_0}$ just above Theorem \ref{thm kato}, and $\Gamma$ in \eqref{claim}, and choose constants $c_1$ and $c_2$ so that $\eta  <c_1 < c_2 < \kappa_{p_0}$

\begin{lemma} \label{lem kato}
There exists $N$ and the value of $\delta_1$ may be chosen sufficiently small that in addition to the conclusions of Theorem \ref{thm kato} we have $|\kappa_z| \ge c_2$ and $\|(T^z_1-T^z_1P^z)^N\|_{\mathbb{C}B_V} \le c_1^N$ for $|z-p_0| \le \delta_1$. 
 \end{lemma}
 
 \n{\bf Proof.}  Choose $c_0$ with $\eta < c_0 < c_1$.  The spectral radius of $T_1^{p_0}-T_1^{p_0}P^{p_0}$ acting on $\mathbb{C}B_V$ is at most $\eta$, so there exists $N$ such that $\|(T_1^{p_0}-T_1^{p_0}P^{p_0})^N\|_{\mathbb{C}V} \le c_0^N$.  The result now follows since $z \to \kappa_z$ and $z \to (T_1^z-T_1^z P^z)^N$ are both analytic functions.
 \qed

\s

Since $z \to \kappa_z$ is analytic and $|\kappa_z| \ge c_2$ for $|z-p_0| < \delta_1$ there is an extension of $\Lambda(p)$ to an analytic function $\Lambda(z) = \log \kappa_z$ for $|z-p_0| < \delta_1$.   Notice $|e^{t \Lambda(z)}| \ge c_2^t$ for $|z-p_0| < \delta_1$.

\s

\begin{lemma}  \label{lem kato1} Given $x \in M$ and $\e > 0$ there exists $\delta_2 \in (0,\delta_1]$ and $t_0 < \infty $such that
 \begin{equation} \label{getz}
  \left| \frac{T_t^z\mathbf{1}(x)}{e^{t \Lambda(z)}} - P^{p_0}\mathbf{1}(x) \right| \le  \e.
  \end{equation}
for all $t \ge t_0$ and $|z-p_0| < \delta_2$.

\end{lemma}

\n{\bf Proof.} 
For $|z-p_0|< \delta_1$ we have 
  $$\|T_{Nk}^z - T_{Nk}^z P^z\|_{\mathbb{C}V} = \|(T_1^z-T_1^zP^z)^{Nk}\|_{\mathbb{C}V} \le c_1^{Nk}.
    $$
Since
   $$
   \left| T_t^z(g_1+ig_2)(x) \right| \le \E^x\left[|g_1(x_t)+ig_2(x_t)| e^{(\Re{z}) A_t}\right]
   $$
then $\|T_t^z\|_{\mathbb{C}V} \le \|T_t^{\Re z}\|_V \le e^{\Gamma t}$ and so    
    $$\|T_{Nk+s}^z - T_{Nk+s}^z P^z\|_{\mathbb{C}V} \le e^{\Gamma s}\|(T_1^z -T_1^zP^z)^{Nk}\|_{\mathbb{C}V} \le e^{\Gamma s}c_1^{Nk}
   \le \left(\frac{e^{\Gamma}}{c_1}\right)^N c_1^{NK+s}
    $$ 
for $0 \le s < N$.  Apply to the function $\mathbf{1}$ at the point $x$:
  \begin{equation} \label{getz2}
  \big|T_t^z\mathbf{1}(x) - T_t^z P^z \mathbf{1}(x) \big| \le \left(\frac{e^{\Gamma}}{c_1}\right)^N c_1^t V(x).
   \end{equation}
Since $T_1^zP^z = \kappa_z P^z = e^{\Lambda(z)}P^z$ then  
   $$
   T_t^z P^z \mathbf{1}(x) = T_s^z \Bigl(T_n^z P^z \mathbf{1}\Bigr)(x) = e^{n \Lambda(z)}T_s^zP^z \mathbf{1}(x) =  e^{t \Lambda(z)}\Bigl(e^{-s\Lambda(z)} T_s^z P^z \mathbf{1}(x)\Bigr)
   $$
where $t = n+s$ with $0 \le s < 1$.  Substituting into \eqref{getz2} and using $|e^{t\Lambda(z)}| \ge c_2^t$ gives
 \begin{equation} \label{getz3}
  \left|\frac{T_t^z\mathbf{1}(x)}{e^{t \Lambda(z)}} - e^{-s\Lambda(z)} T_s^z P^z \mathbf{1}(x) \right| \le \left(\frac{e^{\Gamma}}{c_1}\right)^N \left(\frac{c_1}{c_2}\right)^t V(x)
 \end{equation}  where $t = n+s$ with $0 \le s < 1$,  
   Now consider the function $G(z,s): =e^{-s\Lambda(z)} T_s^z P^z \mathbf{1}(x)$. 
  For $|z-p_0| \le \delta_2 < \delta_1$ and $0 \le s \le 1$ we have 
 \begin{align*}
  \left|P^z\mathbf{1}(x_s) e^{z A_s}\right| & \le \|P^z  \mathbf{1}\|_{\mathbb{C}V} V(x_s) e^{(\Re{z}) A_s}\\
  &  \le \Bigl(\sup_{|z-p_0| \le  \delta_2}\|P^z  \mathbf{1}\|_{\mathbb{C}V}\Bigr) \Bigl( V(x_s) e^{(p_0+\delta_2)A_s}+ V(x_s) e^{(p_0-\delta_2)A_s}\Bigr)
  \end{align*}     
and the right side is  uniformly $\PP^x$-integrable using the same argument as in the proof of Corollary \ref{cor Lam}(iii). 
 Since $(z,s) \mapsto P^z \mathbf{1}(x_s) e^{z A_s}$ is $\PP^x$-almost surely continuous, it follows that
     $(z,s) \mapsto T_s^z P^z \mathbf{1}(x) = \E^x\left[P^z \mathbf{1}(x_s) e^{z A_s}\right]
     $
is a continuous function of $(z,s)$ for $|z-p_0| < \delta_1$ and $0 \le s \le 1$.  Also $z \to \Lambda(z)$ is well-defined and analytic for $|z-p_0| < \delta_1$.  Together $G$ is continuous for $|z-p_0| < \delta_1$ and $0 \le s \le 1$. From Theorem \ref{thm KR} we have $G(p_0,s) = P^{p_0} \mathbf{1}(x) > 0$ for all $s \in[0,1]$. Therefore given $x \in M$ and $\e > 0$ there exists positive $\delta_2 \le \delta_1$ such that 
  $$
 \left|e^{-s\Lambda(z)} T_s^z P^z \mathbf{1}(x) - P^{p_0}\mathbf{1}(x)\right| = \left|G(z,s) - P^{p_0} \mathbf{1}(x)\right| < \frac{\e}{2}
 $$
for $|z-p_0| < \delta_2$ and $0 \le s \le 1$.   Also there exists $t_0$ such that 
  $$
  \left(\frac{e^{\Gamma}}{c_1}\right)^N \left(\frac{c_1}{c_2}\right)^t V(x) < \e/2
  $$
if $t \ge t_0$.  Using these two inequalities in \eqref{getz3} gives \eqref{getz} as required.   \qed

\begin{proposition} \label{prop kato2}  Assume {\bf (H,P,G)} valid for $p = p_0 \pm \delta$ and the same $V$ for some $\delta > 0$.  Given $x \in M$ there exist $\delta_2> 0$ and $t_0$ such that for $t \ge t_0$ the functions $p \to\Lambda(p)$ and $p\to \Lambda_t^x(p) : = \frac{1}{t} \log \E^x\left[e^{pA_t}\right]$ have extensions to analytic functions $\Lambda(z)$ and $\Lambda_t^x(z)$ for $|z-p_0| < \delta_2$.  Moreover there exists a constant $C$ such that
    \begin{equation} \label{Lamt}
      \left|\Lambda_t^x(z) - \Lambda(z)\right| \le \frac{C}{t} 
      \end{equation}
whenever $t \ge t_0$ and $|z-p_0| < \delta_2.$ 
\end{proposition}

\n{\bf Proof.}  In Lemma \ref{lem kato1} choose $\e = P^{p_0}\mathbf{1}(x)/2$.  Then \eqref{getz} implies
    $$
   \left| \frac{T_t^z\mathbf{1}(x)}{e^{t\Lambda(x)}}\right| \ge P^{p_0}\mathbf{1}(x)/2 > 0
   $$
for $t \ge t_0$ and $|z-p_0| < \delta_2$. 
For fixed $x$ and $t \ge t_0$ the function $z \to T_t^z\mathbf{1}(x)$ is analytic and bounded away from 0, so that $z  \to \Lambda_t^x(z):=t^{-1} \log T_t^z \mathbf{1}(x)$ is well-defined and analytic for $|z-p_0| < \delta_2$.  
Also  \eqref{getz} implies
  $$
  \left|t \Lambda_t^x(z) - t \Lambda(z)\right| = \left| \log \left(\frac{T_t^z \mathbf{1}(x)}{e^{ t\Lambda(z) }}\right) \right| \le \sup\{|\log w|: |w-P^{p_0}\mathbf{1}(x)| \le P^{p_0}\mathbf{1}(x)/2\}
  $$
whenever $t \ge t_0$ and $|z-p_0| < \delta_2.$   This gives \eqref{Lamt} with $C = \sup\{|\log w|: |w-P^{p_0}\mathbf{1}(x)| \le P^{p_0}\mathbf{1}(x)/2\} < \infty.$  \qed

\s

\n {\bf Proof of Theorem \ref{thm kato2}.}  For $k = 0$ the inequality \ref{deriv} is given by \eqref{Lamt}in Proposition \ref{prop kato2} for $|z-p_0| < \delta_2$, and then the inequalities for $k \ge 1$ on a smaller disc $|z-p_0| < \delta_3$ for $0 < \delta_3 < \delta_2$ follow from the Cauchy integral formula for the analytic function $z \mapsto \Lambda_t^x(z) - \Lambda(z)$.  \qed

\section{Proofs for Section \ref{sec near}} \label{sec near proof}

\n{\bf Proof of Lemma \ref{lem nu}.}  Combining the growth assumptions {\bf (1)} at $p = \pm \delta$ implies
  $$
  \frac{LV(x)}{V(x)} + \delta\left|\frac{YV(x)}{V(x)} + Q(x) \right| + \frac{\delta^2}{2}R(x) \to -\infty
  $$
as $x \to \infty$.  Let $U = \log V$, then
  $$
  LU(x) + \frac{1}{2} \sum_{j=1}^m \bigl(X_jU(x)\bigr)^2+ \delta\left|YU(x) + Q(x) \right| + \frac{\delta^2}{2}R(x) \to -\infty
  $$ 
as $x \to \infty$.  By Meyn and Tweedie \cite[Theorem 4.2]{MTIII} we deduce that 
$\frac{1}{2} \sum_{j=1}^m \bigl(X_jU(x)\bigr)^2+ \delta\left|YU(x) + Q(x) \right| + \frac{\delta^2}{2}R(x)$ is $\nu$-integrable.  It follows that $X_jU \in L^2(\nu)$ and $q_j \in L^2(\nu)$ for $1 \le j\le m$, and $YU+Q \in L^1(\nu)$.  Also
  $$
  |YU(x)| = \left|\sum_{j=1}^m q_j(x)X_jU(x) \right| \le \frac{1}{2}\sum_{j=1}^m \left((q_j(x)^2+(X_jU(x))^2 \right).
  $$
Therefore $YU \in L^1(\nu)$ and then $Q = (YU+Q) - YU \in L^1(\nu)$. 
\qed

\subsection{``Abstract'' method} \label{sec near abstract}

\n{\bf Proof of Theorem \ref{thm clt}(i).}  By Theorem \ref{thm kato2} for $|p| < \delta$ we have
   $$
   \frac{1}{t} \frac{\E^x[A_t e^{pA_t}]}{\E^x[e^{pA_t}]} = \big(\Lambda_t^x\bigr)'(p) \to \Lambda'(p)
   $$
and 
  $$
   \frac{1}{t}\left( \frac{\E^x[A_t^2 e^{pA_t}]}{\E^x[e^{pA_t}]}-  \frac{\bigl(\E^x[A_t e^{pA_t}]\bigr)^2}{\bigl(\E^x[e^{pA_t}]\bigr)^2}\right) = \big(\Lambda_t^x\bigr)''(p) \to \Lambda''(p)
   $$
as $t \to \infty$.  Putting $p = 0$ gives
   $$
   \frac{1}{t} \E^x[A_t] \to \Lambda'(0) \quad \mbox{ and } \quad \frac{1}{t} \mbox{Var}(A_t) \to \Lambda''(0)
   $$   
as $t \to \infty$.  The second limit implies $\frac{1}{t}A_t - \frac{1}{t}\E^x[A_t] \to 0$ in $\PP^x$-probability as $t \to \infty$, and then the first limit implies $\frac{1}{t}A_t \to \Lambda'(0)$ in $\PP^x$-probability as $t \to \infty$.   Comparing with \eqref{lln} we deduce $\Lambda'(0) = \lambda$.  \qed

\s

\n {\bf Proof of Theorem \ref{thm clt}(ii)}  The method is based closely on results of Li Ming Wu \cite{WuL95}.  Theorem \ref{thm kato2} ensures that the distribution of $A_t/t$ satisfies Wu's condition of $C^2$-regularity, and then the following result is due to Li Ming Wu \cite[Theorem 1.2]{WuL95}.  We have changed the centering from $A_t - \E^x[A_t]$ to $A_t - t \Lambda'(0)$ and adapted the statement and proof to our notation.

\begin{proposition}  \label{prop wu} Assume {\bf (H,P,G)} valid for two values $p = \pm \delta$ with same function $V$ for some $\delta > 0$.  Suppose either (i) $b_t = \sqrt{t}$ or else (ii) $b_t/\sqrt{t} \to \infty$ and $b_t/t \to 0$.  In both cases write $a_t = b_t^2/t$.  Then 
   \begin{equation} \label{wu}
  \lim_{t \to \infty} \frac{1}{a_t}\log \E^x\left[\exp\left\{\frac{pa_t (A_t- t\Lambda'(0))}{b_t}\right\}\right] 
  = \frac{1}{2}p^2 \Lambda''(0)
   \end{equation}
fo all $p\in \R$.    
\end{proposition}

\n{\bf Proof.}  For fixed $p \in \R$ we have  
  \begin{align*}
   \frac{1}{a_t}\log \E^x\left[\exp\left\{\frac{p a_t(A_t-t \Lambda'(0))}{b_t}\right\}\right]
    & =  \frac{1}{a_t}\left(\log \E^x\left[\exp\left\{\frac{p a_t  A_t}{b_t}\right\}\right] - \frac{ p a_t t\Lambda'(0)}{b_t}\right) \\
    & = \frac{t}{a_t}\left(\Lambda_t^x\bigl(\frac{p a_t}{b_t}\bigr) - \frac{p a_t}{b_t}\Lambda'(0)\right)\\
          &  = \frac{t^2}{b_t^2}\left(\Lambda_t^x\bigl(\frac{pb_t}{t}\bigr) - \frac{p b_t}{t}\Lambda'(0)\right)\\
        &  = \frac{t^2}{b_t^2}\left(\Lambda_t^x\bigl(\frac{pb_t}{t}\bigr) - \frac{p b_t}{t}(\Lambda_t^x)'(0)\right)
              + \frac{p t}{b_t}\bigl((\Lambda_t^x)'(0) - \Lambda'(0)\bigr).
    \end{align*}
The second term on the right tends to 0 as $t \to \infty$ by \eqref{deriv} in Theorem \ref{thm kato2}.  The first term on the right is 
   \begin{align*}
  \frac{t^2}{b_t^2}\left(\Lambda_t^x\bigl(\frac{pb_t}{t}\bigr) - \frac{p b_t}{t}(\Lambda_t^x)'(0)\right)
    & = \frac{t^2}{b_t^2} \int_0^{p b_t/t} (\Lambda_t^x)''(s) \bigl(\frac{pb_t}{t} - s\bigr)ds\\
     & = \int_0^p (\Lambda_t^x)''\bigl(\frac{b_tu}{t}\bigr) (p -u )du\\
     & \to  \int_0^p \Lambda''(0) \bigl(p -u \bigr)du = \frac{1}{2}p^2 \Lambda''(0)
    \end{align*}
where the first equality is valid for sufficiently small $|p b_t/t|$, and the convergence again uses \eqref{deriv} in Theorem \ref{thm kato2}.  \qed

\s

 Taking $b_t = \sqrt{t}$ gives 
$$
  \lim_{t \to\infty}  \E^x\left[\exp\left\{p\frac{(A_t - t \Lambda'(0))}{\sqrt{t}}\right\}\right]
     = \exp\left\{\frac{1}{2}p^2 \Lambda''(0)\right\} = \E\left[e^{p Z}\right] \quad \mbox{ for all } p \in \R,
     $$
where $Z \sim N(0,\Lambda''(0))$. It follows, see for example \cite[Lemma 1.3(b)]{WuL95}, that $
   \dfrac{(A_t-t \Lambda'(0))}{\sqrt{t}}$ converges in distribution to $N(0,\Lambda''(0))$.   This completes the proof of Theorem \ref{thm clt}(ii).  \qed 
  
\s

\n{\bf Proof of Theorem \ref{thm mod dev}.}
Recall the G\"{a}rtner-Ellis theorem:  if
    $$
  \lim_{t \to\infty}  \frac{1}{t} \log \E^x[e^{pA_t}]  = \Lambda(p) < \infty
    $$
for all $p \in \R$ then $A_t/t$ satisfies the large deviations principle with speed $t$ and rate function $I(s)= \sup_{p \in \R} (ps-\widetilde{\Lambda}(p))$.  See for example Dembo and Zeitouni \cite{DZ10}.  Now rescale in time and space using $a_t \to \infty$ and $b_t \to \infty$.  If
   $$
    \frac{1}{a_t} \log \E^x[e^{a_t p(A_t-t \Lambda'(0))/b_t}] \to \widehat{\Lambda}(p)
    $$
then the G\"{a}rtner-Ellis theorem implies $(A_t-t\Lambda'(0))/b_t$ satisfies the large deviations principle with speed $a_t$ and rate function $\widehat{I}(x) = \sup_{p \in \R} (px-\widehat{\Lambda}(p))$.  From Proposition \ref{prop wu} we have $a_t = b_t^2/t \to \infty$ and $\widehat{\Lambda}(p) = p^2\Lambda''(0)/2$.  By assumption $\Lambda''(0) \neq 0$ and so $\widehat{I}(s) = \sup_{p \in \R}(ps - p^2 \Lambda''(0)/2) = s^2/(2\Lambda''(0))$.   \qed

\subsection{Using eigenfunctions} \label{sec near eigen}

Recall the eigenvalue/eigenfunction characterization of $\Lambda(p)$ in Corollary \ref{cor Lam}:
\begin{equation} \label{char}
     \left(L + pY + pQ + \frac{p^2}{2}R\right)\phi_p(x) = \Lambda(p)\phi_p(x)
     \end{equation}
where $Y(x) =\sum_{j=1}^m q_j(x)X_j(x)$ and $R(x) =  \sum_{j=1}^m \bigl(q_j(x)\bigr)^2$.  By Theorem \ref{thm kato1} the mappings $p \to \Lambda(p)$ and $p \to \phi_p$ are real analytic.   
     Differentiating with respect to $p$ gives
   \begin{equation} \label{char1}
     \left(Y + Q + pR\right)\phi_p(x) + \left(L + pY + pQ + \frac{p^2}{2}R\right)\frac{\partial \phi_p}{\partial p}(x) = \Lambda'(p)\phi_p(x) + \Lambda(p)\frac{\partial \phi_p}{\partial p}(x), 
    \end{equation}
and differentiating again gives
   \begin{align}
  \lefteqn{ R(x)\phi_p(x) + 2(Y+Q+pR)\frac{\partial \phi_p}{\partial p}(x) + \left(L + pY + pQ + \frac{p^2}R\right)\frac{\partial^2 \phi_p}{\partial p^2}(x) } \hspace{25ex} \nonumber \\
  & = \Lambda''(p)\phi_p(x) + 2\Lambda'(p)\frac{\partial \phi_p}{\partial p}(x) +  \Lambda(p)\frac{\partial^2 \phi_p}{\partial p^2}(x). \label{char2} 
  \end{align}     
Recall from Corollary \ref{cor Lam}(ii) that $\phi_p \in C^\infty(M)$.  Using the hypoellipticity of $L_p$ and \eqref{char1} we have $\frac{\partial \phi_p}{\partial p} \in C^\infty(M)$.  Repeating the argument on \eqref{char2}, we have $\frac{\partial^2 \phi_p}{\partial p^2} \in C^\infty(M)$.

\s

Write $\Phi_1= \frac{\partial \phi_p}{\partial p}\big|_{p=0}$ and $\Phi_2(x)= \frac{\partial^2 \phi_p}{\partial p^2}\big|_{p=0}$, and note that $\Phi_1$ and $\Phi_2$ are both in $B_V$.  
 Putting $p = 0$ in \eqref{char} gives $L\phi_0(x) = \Lambda(0)\phi_0(x)$ so that $\Lambda(0) = 0$ and we can take $\phi_0(x) \equiv 1$.  Then putting $p= 0$ in \eqref{char1} gives
  \begin{equation} \label{poisson}
 Q(x) + L\Phi_1(x) = \Lambda'(0).
 \end{equation} 
This is the so-called Poisson equation, and the existence of the solution $\Phi_1$ is crucial in obtaining a central limit theorem.  See for example \cite{Bha82} and \cite{GM96}.  Putting $p=0$ in \eqref{char2} gives
   $$
   R(x) + 2(Y+Q)\Phi_1(x) + L\Phi_2(x) 
    = \Lambda''(0) + 2\Lambda'(0)\Phi_1(x)
    $$
  Then $\Phi_1$ and $\Phi_2$ are in $B_V\cap C^\infty(M)  $ and satisfy
     \begin{align*}
     L\Phi_1(x) & = \Lambda'(0) - Q(x) \\
     L\Phi_2(x) & = \Lambda''(0) +2\bigl(\Lambda'(0)-Q(x)\bigr)\Phi_1(x) - 2Y\Phi_1(x) - R(x).
     \end{align*}    
Using It\^{o}'s formula we have
\begin{align} \nonumber
   A_t+ \Phi_1(x_t) & = \Phi_1(x_0) + \int_0^t \bigl(Q(x_s) + L\Phi_1(x_s)\bigr)ds + \sum_{j=1}^m  \int_0^t \bigl(q_j(x_s)+ X_j\Phi_1(x_s)\bigr)dW_s^j \\
       & =  \Phi_1(x_0) + \Lambda'(0)t + M_t \label{pin}
  \end{align}
  where 
  $$
   M_t = \sum_{j=1}^m \int_0^t   \bigl(q_j(x_s)+ X_j \Phi_1(x_s)\bigr)dW_s^j
   $$
is a local martingale.  
    An early version of \eqref{pin} for a linear SDE is due to Pinsky \cite[Eqn (9)]{pin74}.  
To proceed we need more information about the function $\Phi_1$ and the quadratic variation $\langle M \rangle_t$ 

\begin{lemma} \label{lem clt}   Assume {\bf (H,P,G)} valid for two values $p = \pm \delta$ with same function $V$ for some $\delta > 0$.  Then
   
   (i) $(\Phi_i)^q \in B_V \subset L^1(\nu)$ for all $q \ge 1$ and $i = 1,2$;
   
   (ii) $$
   \int_M Q d \nu = \Lambda'(0); 
   $$
   
   (iii) 
     $$
     \sum_{j=1}^m \int_M \bigl(q_j(x)+X_j \Phi_1(x)\bigr)^2d\nu(x) = \Lambda''(0).
     $$
     \end{lemma}

\n{\bf Proof.} Writing $U = \log V$ and letting $0 < \alpha <1$ gives 
  \begin{align*}
   \frac{L_{\alpha p}V^\alpha(x)}{V^\alpha(x)} 
    & = \alpha LU(x)+ \frac{\alpha^2}{2} \sum_{j=1}^r \bigl(X_jU(x)\bigr)^2+ \alpha^2 p \sum_{j=1}^r q_j(x)X_jU(x) + \alpha pQ(x) + \frac{\alpha^2p^2}{2}\sum_{j=1}^r q_j^2(x)\\
    & = \alpha LU(x)+ \frac{\alpha^2}{2} \sum_{j=1}^r \bigl(X_jU(x)+ p q_j(x) \bigr)^2+ \alpha pQ(x) \\
    & \le \alpha \left(LU(x)+ \frac{1}{2} \sum_{j=1}^r \bigl(X_jU(x)+ p q_j(x) \bigr)^2+ pQ(x) \right)\\
    & = \alpha\left(\frac{L_pV(x)}{V(x)}\right).
    \end{align*}  
Since {\bf (G)} is satisfied for $p = \pm \delta$ with the function $V$, then it is also satisfied for  $p = \pm \alpha\delta$ with the function $V^\alpha.$  
Since $\Phi_1$ and $\Phi_2$ do not depend on the value of $\delta$, we obtain $\Phi_1 \in B_{V^\alpha}$ and $\Phi_2 \in B_{V^\alpha}$.  
 Taking $\alpha = 1/q$ gives $(\Phi_i)^q \in B_V \subset L^1(\nu)$ for $i = 1,2$, and (i) is proved.   
 
Next, we have $\Phi_1 \in L^1(\nu)$ and $L\Phi_1 = \Lambda'(0) - Q  \in L^1(\nu)$ and $\Phi_1(x) /V(x) \to 0$  as $x \to \infty$.  By Baxendale \cite[Proposition B.1]{bax24} we deduce 
  $$
  \int_M \Bigl(\Lambda'(0) - Q)\Bigr)d\nu = \int_M L \Phi_1 d\nu =0.
  $$ 


Finally
   \begin{align}
   L(\Phi_2-\Phi_1^2)(x) & = \Lambda''(0) +2\bigl(\Lambda'(0)-Q(x)\bigr)\Phi_1(x) - 2Y\Phi_1(x) - R(x) \nonumber \\
   & \quad  - 2\Phi_1(x) L\Phi_1(x) - \sum_{j=1}^m \bigl(X_j\Phi_1(x)\bigr)^2 \nonumber\\
   & = \Lambda''(0)- 2Y\Phi_1(x) - R(x) - \sum_{j=1}^m \bigl(X_j\Phi_1(x)\bigr)^2 \nonumber\\
   & = \Lambda''(0)- \sum_{j=1}^m\bigl(q_j(x)+X_j \Phi_1(x) \bigr)^2. \label{LF}
    \end{align} 
At this stage we would like to integrate \eqref{LF} with respect to $\nu$ and claim $\int_M L(\Phi_2-\Phi_1^2)d\nu = 0$.  From (i) above we have $\Phi_2-\Phi_1^2 \in L^1(\nu)$ and $|\Phi_2-\Phi_1^2|/V(x) \to 0$ as $x \to \infty$.   Therefore we may complete the proof of Lemma \ref{lem clt} using the following variation on Baxendale \cite[Proposition B.1]{bax24}.  \qed

\begin{proposition} \label{prop BG}
  Let $\{x_t:\ t\geq 0\}$ be a non-explosive 
diffusion process on a $\sigma$-compact manifold $M$ with
invariant probability measure $\nu$. Let ${\cal L}$ be an operator acting
on $C^2(M)$ functions that agrees with the generator of $\{x_t:\
t\geq 0\}$ on $C^2$ functions with compact support.  Let $F \in
C^2(M)$.  Suppose $F$ is $\nu$-integrable and ${\cal L}F$ is bounded below.  Suppose also there exists a positive $V \in C^2(M)$ satisfying
${\cal L} V(x) \leq kV(x)$ for some $k<\infty$ such that $|F(x)|/V(x) \to 0$ as $|F(x)| \to \infty$. Then ${\cal L}F$ is $\nu$-integrable and 
 $$
 \int_M {\cal L}F(x)\, d\nu(x) = 0.
 $$
 \end{proposition}

 \n{\bf Proof.} 
  We adopt the notation and much of the proof of \cite[Proposition B.1]{bax24}.  The fact that ${\cal L}F$ is bounded below implies that $\E^x \int_0^{t \wedge \delta_n} {\cal L}F(X_s)ds \to \E^x \int_0^t {\cal L}F(X_s)ds \le \infty$ as $n \to \infty$, using a combination of bounded and monotone convergence theorems.  This gives
      $$
 \E^x(F(X_t)) - F(x) =
   \E^x \left( \int_0^t {\cal L}F(X_s) ds \right) \le \infty
 $$ 
for all $x \in M$.  Integrating and using the invariance of $\nu$, together with Fubini's theorem on the right side gives
$$
 0 = \int_M  F(x) \,d\nu(x) -  \int_M F(x) \,d\nu(x) = t\int_M {\cal L}F(x)\,d\nu(x).
  $$
In particular $\int_M {\cal L}F(x)\,d\nu(x) <\infty$, so that ${\cal L}F$ is $\nu$-integrable,  and the proof is completed. \qed

 \s
\n{\bf Proof of Theorem \ref{thm clt}}   The result (i) follows from Lemma \ref{lem clt} together with \eqref{lln}.  For (ii) we adapt the arguments in \cite[Prop 2.6]{pin74} and \cite[Cor 3.2]{Bax87} to the non-compact setting.  Recall the local martingale $M_t$ in \eqref{pin}.
We have
  $$
      \frac{1}{t}\langle M\rangle _t = \sum_{j=1}^m \frac{1}{t}\int_0^t \bigl(q_j(x_s) + X_j \Phi_1(x_s)\big)^2\,ds \to \sum_{j=1}^m  \int_M\bigl(q_j(x)+X_j\Phi_1(x)\bigr)^2d\nu(x)  = \Lambda''(0)<\infty.
      $$
 If $\Lambda''(0) > 0$ then $\langle M\rangle _t \to \infty$ almost surely and we may write $M_t = B_{\langle M\rangle_t}$ where $B$ is a standard Brownian motion.  See for example Revuz and Yor \cite[Theorem V.1.6]{RevYor}.  Then \eqref{pin} gives 
    $$
    \sqrt{t}\left(\frac{A_t}{t} - \Lambda'(0)\right) = \frac{\Phi_1(x_0)}{\sqrt{t}} - \frac{\Phi_1(x_t)}{\sqrt{t}} + \frac{1}{\sqrt{t}}B_{\langle M \rangle_t}. 
    $$
Since $\Phi_1 \in B_V$ then $\E |\Phi_1(x_t)|$ is bounded as $t \to \infty$ and so 
     $$
    \frac{\Phi_1(x_0)}{\sqrt{t}} - \frac{\Phi_1(x_t)}{\sqrt{t}} \to 0
     $$
in $\PP^x$-probability as $t \to \infty$.  The result now follows since
     $$
     \frac{1}{\sqrt{t}} B_{\langle M \rangle_t} = \sqrt{\frac{\langle M \rangle_t}{t} } \frac{B_{\langle M \rangle_t}}{\sqrt{\langle M \rangle_t}}  \to \sqrt{\Lambda''(0)} N\left(0,1\right) = N(0,\Lambda''(0)) 
     $$
in distribution as $t \to \infty$. If $\Lambda''(0) = 0$ then $M_t = 0$ for all $t \ge 0$ and so
   $$
    \sqrt{t}\left(\frac{A_t}{t} - \Lambda'(0)\right) = \frac{\Phi_1(x_0)}{\sqrt{t}} - \frac{\Phi_1(x_t)}{\sqrt{t}} \to 0 
    $$
in $\PP^x$-probability as $t \to \infty$. \qed     

\s

\n{\bf Proof of Theorem \ref{thm linear}.} (i) implies (ii) by Lemma \ref{lem nu}, and (ii) implies (iii) by \eqref{pin}.  Given (iii) we have 
  $$
   T_t^p (e^{p \Phi_1})(x)  = \E^x\left[ e^{p\Phi_1(x_t)} e^{pA_t}\right] 
   = \E^x\left[e^{p(\lambda t + \Phi(x_0))}\right] = e^{p \lambda t} e^{p \Phi_1(x)}.
 $$
This suggests $e^{p \Phi_1}$ is an eigenfunction for $T_1^p$ corresponding to the eigenvalue $e^{p\lambda}$, and hence that $\Lambda(p) = p\lambda$.  But we need to check that $e^{p \Phi_1}$ is in the Banach space $B_V$. For $|p| <\delta$ we know that $T_t^p$ acting on $B_V$ has an eigenvalue $\kappa_p^t$ and corresponding positive eigenfunction $\phi_p \in B_V$.  Fix a non-empty open set $U \subset M$ with compact closure and write $g = \mathbf{1}_U$.  Fix $x \in M$.  By Theorem \ref{thm KR} we have $T_t^p g(x) \sim \kappa_p^t \phi_p(x) \widehat{\nu}_p(U)$ as $t \to \infty$ with $\phi_p(x) \widehat{\nu}_p(U) > 0$.  Also, there exists $c$ such that $g \le c e^{p \Phi_1}$ and so 
   $$
   T_t^p g(x) \le c T_t^p(e^{p \Phi_1})(x) = ce^{p \lambda t} e^{p \Phi_1(x)}.
   $$
Together this gives $\kappa_p \le e^{p \lambda}$, and hence $\Lambda(p) \le p \lambda$.  Since $\Lambda(p)$ is convex on $(-\delta,\delta)$ and $\Lambda'(0) = \lambda$, we deduce $\Lambda(p) = \lambda p$ for all $p \in (-\delta,\delta)$, giving (iv). Clearly (iv) implies (i), and the proof is complete.  \qed

\section{Proofs for Section \ref{sec asymp}} \label{sec asymp proof}

In each case we obtain an upper bound and a lower bound separately.

 \subsection{Proof of Proposition \ref{prop asymp}(i): $Q(x) = x^2$} \label{Q2}

For the upper bound, take $V(x) = e^{\gamma x^2}$.  Then 
 \begin{align*}
   \frac{L_pV(x)}{V(x)} 
    & = -2 \gamma b x^4 +(2 \gamma a +2\gamma^2 \sigma^2 )x^2 + \gamma \sigma^2 + px^2\\
    & = -2b \gamma x^4+ (2\gamma a +2 \gamma^2 \sigma^2 +p)x^2+ \gamma \sigma^2.
    \end{align*}
We get {\bf (G)} with any $\gamma > 0$ and this will work for all $p$.  We have 
  \begin{align*}
   \Lambda(p) \le \sup_{x \in \R} \frac{L_pV(x)}{V(x)} & = \sup_{x \in \R}\Bigl(-2b \gamma x^4+ (2\gamma a +2 \gamma^2 \sigma^2 +p)x^2+ \gamma \sigma^2\Bigr) \\
    & = \frac{(2\gamma a +2 \gamma^2 \sigma^2 +p)^2}{8b \gamma} + \gamma \sigma^2.
  \end{align*}
where the last equality is true so long as $2 \gamma a +2 \gamma^2 \sigma^2 +p \ge 0$.  In particular it will be true for sufficiently large $\gamma$. 
Recall we can use any value of $\gamma> 0$ with any value of $p >0$.  Putting $\gamma = c \sqrt{p}$ with $c> 0$ we get
     $$
  \limsup_{p \to\infty}\frac{\Lambda(p)}{p^{3/2}}  \le  \frac{(2 c^2 \sigma^2 +1)^2}{8b c}.
  $$
The right side is minimized when $6c^2 \sigma^2 = 1$ and this gives
   $$
    \limsup_{p \to\infty}\frac{\Lambda(p)}{p^{3/2}}  \le  \left(\frac{2}{3}\right)^{3/2} \frac{\sigma}{b}.
    $$

For the lower bound, consider $g(x) = (x^2)^A$ for some $A > 1$.  Then 
  \begin{align*}
    (L_pg(x) 
      & = \Bigl(2Aa - 2Abx^2 +\sigma^2 A(2A-1)x^{-2} + px^2 \Bigr)g(x)\\
      & \ge \widetilde{\Gamma}_{p,A}g(x)
      \end{align*}       
where
   \begin{align*}
   \widetilde{\Gamma}_{p,A} & = \inf_{x \in\R}\Bigl(2Aa - 2Abx^2 +\sigma^2 A(2A-1)x^{-2} + px^2 \Bigr)\\
     & = 2Aa + \inf_{x \in\R}\Bigl((p- 2Ab)x^2 +\sigma^2 A(2A-1)x^{-2} \Bigr)\\
     & = 2Aa + 2\sqrt{p-2Ab} \sqrt{\sigma^2 A(2A-1)}
  \end{align*}
so long as $2Ab < p$.  For large $p$ choose $A = p/(3b)$.  Then 
  $$
  \widetilde{\Gamma}_{p,A} \sim \left(\frac{2p}{3}\right)^{3/2} \frac{\sigma}{b}
  $$
as $p \to \infty$, and so 
  $$
   \liminf_{p \to \infty} \frac{\Lambda(p)}{p^{3/2}} \ge \liminf_{p \to \infty} \frac{A_{p,A}}{p^{3/2}} \ge \left(\frac{2}{3}\right)^{3/2} \frac{\sigma}{b}.
    $$
The upper and lower asymptotic estimates coincide, and the proof is complete. \qed

 \subsection{Proof of Proposition \ref{prop asymp}(ii): $Q(x) = x^4$} \label{Q4}

For the upper bound, again take $V(x) = e^{\gamma x^2}$.  Then 
 \begin{align*}
   \frac{LpV(x)}{V(x)} 
       & = (p -2b \gamma) x^4+ (2\gamma a +2 \gamma^2 \sigma^2 )x^2+ \gamma \sigma^2.
    \end{align*}
For each $p$ the condition {\bf (G)} is satisfied so long as $\gamma > p/(2b)$.  For $\gamma > p/(2b)$ we have
   \begin{align*}
   \Lambda(p) \le \sup_{x \in \R} \frac{L_p V(x)}{V(x)} & = \sup_{x \in \R}\Bigl((p -2b \gamma) x^4+ (2\gamma a +2 \gamma^2 \sigma^2 )x^2+ \gamma \sigma^2\Bigr) \\
     & = \frac{(2\gamma a +2 \gamma^2 \sigma^2 )^2}{4(2b\gamma - p)} + \gamma \sigma^2,
     \end{align*} 
where the last equality is true so long as $2 \gamma a +2 \gamma^2 \sigma^2 \ge 0$.  In particular it will be true for sufficiently large $\gamma$. Choose $\gamma = cp$ with $c > 1/(2b)$.  Then 
   $$
    \limsup_{p \to \infty} \frac{\Lambda(p)}{p^3} \le \frac{\sigma^4 c^4}{(2bc-1)}
    $$
 The right side is minimized when $c = 2/(3b)$ and then
   $$
   \limsup_{p \to \infty} \frac{\Lambda(p)}{p^3} \le \frac{16}{27}\left(\frac{\sigma}{b}\right)^4.
   $$

For the lower bound, again consider $g(x) = (x^2)^A$ for some $A > 1$.  Then 
      \begin{align*}
    L_pg(x)
           & = \Bigl(2Aa - 2Abx^2 +\sigma^2 A(2A-1)x^{-2} + px^4 \Bigr)g(x)\\
      & \ge \widetilde{\Gamma}_{p,A}g(x)
      \end{align*}       
where
 $$
  \widetilde{\Gamma}_{p,A}  = \inf_{x \in\R}\Bigl(2Aa - 2Abx^2 +\sigma^2 A(2A-1)x^{-2} + px^4 \Bigr),
  $$
For large $p$ choose $A = cp^2$ for some $c > 0$.  Then putting $y=x^2/p$ we get
    \begin{align*}
   \widetilde{\Gamma}_{p,A} & = \inf_{y >0}\Bigl(2cp^2a - 2cp^3by +\sigma^2 cp(2cp^2-1)y^{-1} + p^3y^2)\\
     & \sim p^3 \inf_{y > 0}\Bigl( - 2cb y +2\sigma^2 c^2y^{-1} + y^2 \Bigr)
  \end{align*}
as $p \to \infty$.  Next choose $c$ so as to maximize the right side above.  An elementary saddle-point analysis gives $c = 8 \sigma^2/(9 b ^3)$.  With this value of $c$ and then the substitution $y = 4 \sigma^2 z/(3b^2)$  we have
  \begin{align*}
  \liminf_{p \to \infty} \frac{\Lambda(p)}{p^3} \ge \liminf_{p \to \infty}\frac{\widetilde{\Gamma}_{p,A}}{p^3}
  & = \inf_{y > 0}\Bigl( - 2cb y +2\sigma^2 c^2y^{-1} + y^2 \Bigr) \\
  & = \frac{16 \sigma^4}{27 b^4} \inf_{z > 0}\Bigl( -4z+2z^{-1} + 3z^2 \Bigr) = \frac{16 \sigma^4}{27 b^4}.
  \end{align*}
The upper and lower asymptotic estimates coincide, and the proof is complete. \qed  

\subsection{Proof of Remark \ref{rem corr}}  \label{sec corr}  With $q(x) = x$ and correlation $\rho \in (-1,1]$ we have $Y(x) = \rho \sigma x \frac{d}{dx}$ and $R(x) = x^2$.  For the upper bound, take $V(x) = e^{\gamma x^2}$.  Then  
 \begin{align*}
  \frac{L_pV(v)}{V(x)} 
    &=  -2b \gamma x^4+ (2\gamma a +2 \gamma^2 \sigma^2 )x^2+ \gamma \sigma^2 + 2 \rho \gamma p \sigma x^2 + \frac{1}{2}p^2x^2 \\
       & = -2b \gamma x^4+ (2\gamma a +2 \gamma^2 \sigma^2 +2 \rho\gamma p \sigma+\frac{1}{2}p^2)x^2+ \gamma \sigma^2.
    \end{align*}
We get {\bf (G)} with any $\gamma > 0$ and this will work for all $p$. We have 
 \begin{align*}
  \Lambda(p) \le \sup_{x \in M} \frac{L_pV(x)}{V(x)} & = \sup_{x \in \R}\Bigl( -2b \gamma x^4+ (2\gamma a +2 \gamma^2 \sigma^2+2 \rho \gamma p \sigma +\frac{1}{2}p^2)x^2+ \gamma \sigma^2\Bigr)\\
  & = \frac{(2\gamma a +2 \gamma^2 \sigma^2+2 \rho \gamma p \sigma +\frac{1}{2}p^2)^2}{8b \gamma} + \gamma \sigma^2,
  \end{align*}
where the last equality is true so long as $2\gamma a +2 \gamma^2 \sigma^2+2 \rho \gamma p \sigma +\frac{1}{2}p^2 \ge 0$.  In particular since $\rho > -1$ it will be true for sufficiently large $\gamma$.  Recall we can use any value of $\gamma> 0$ with any value of $p >0$.  Putting $\gamma = c p/(2\sigma)$ for some $c  >0$ we get
     $$
  \limsup_{p \to\infty}\frac{\Lambda(p)}{p^3}  \le  \frac{\sigma}{16 b}  \frac{(c^2+2 \rho c+1)^2}{c}.
  $$
Choosing $c = (\sqrt{3+\rho^2} -\rho)/3$ to minimize the right side gives
  \begin{equation} \label{upper}
  \limsup_{p \to\infty}\frac{\Lambda(p)}{p^3}  \le 
  \frac{\sigma}{27b}(\sqrt{3+\rho^2}-\rho)(\sqrt{3+\rho^2}+2 \rho)^2.
  \end{equation}
For the lower bound, again consider $g(x) = (x^2)^A$ for some $A > 1$.  Then 
   \begin{align*}
    L_pg(x) & = Lg(x)+pYg(x) + \frac{1}{2}p^2x^2g(x) \\
    & = \Bigl(2Aa - 2Abx^2 +\sigma^2 A(2A-1)x^{-2} +2A \rho p\sigma + \frac{1}{2}p^2x^2 \Bigr)g(x)\\
         & = \Bigl(2A(a+ \rho p\sigma) +(\frac{1}{2}p^2-2Ab)x^2 +\sigma^2 A(2A-1)x^{-2} \Bigr)g(x)\\
      & \ge \widetilde{\Gamma}_{p,A}g(x)
      \end{align*}       
where
   \begin{align*}
   \widetilde{\Gamma}_{p,A} & = \inf_{x \in\R}\Bigl(2A(a+ p\sigma) - 2Abx^2 +\sigma^2 A(2A-1)x^{-2} + \frac{1}{2}p^2x^2 \Bigr)\\
     & = 2A(a+ \rho p\sigma) + \inf_{x \in\R}\Bigl((\frac{1}{2}p^2- 2Ab)x^2 +\sigma^2 A(2A-1)x^{-2} \Bigr)\\
     & = 2A(a+\rho p\sigma) + 2\sqrt{p^2/2-2Ab} \sqrt{\sigma^2 A(2A-1)}
       \end{align*}
so long as $4Ab < p^2$.  For large $p$ choose $A = cp^2/(4b)$ for some $c$ with $0 < c < 1$.  Then 
 $$
  \liminf _{p \to \infty}\frac{\Lambda(p)}{p^3} \ge \liminf_{p \to \infty}\frac{ \widetilde{\Gamma}_{p,A}}{p} = \frac{\sigma}{2b}\Bigl(\rho c + c\sqrt{1-c}\Bigr).
     $$ 
Write $c = 1-\widetilde{c}^2$ and then choose $\widetilde{c} = (\sqrt{3+\rho^2} - \rho)/3$ to maximize the right side.  Then 
   \begin{equation} \label{lower}
   \liminf_{p \to \infty} \frac{\Lambda(p)}{p^3} \ge  \frac{\sigma}{27b}(\sqrt{3+\rho^2}-\rho)(\sqrt{3+\rho^2}+2 \rho)^2.
   \end{equation} 
The upper and lower asymptotic estimates coincide, and the proof is complete. \qed


\subsection{Proof of Proposition \ref{prop asymp}(iii,iv)}  

Proposition \ref{prop asymp}(iii) corresponds to correlation coefficient $\rho = 1$ in Remark \ref{rem corr} and this case is included in Section \ref{sec corr} above.  It remains to see how the estimates change in the special case $\rho = -1$.  For the upper bound, take $V(x) = e^{\gamma x^2}$.  With $\rho = -1$ we have  
 \begin{align*}
  \frac{L_pV(v)}{V(x)} 
    &=  -2b \gamma x^4+ (2\gamma a +2 \gamma^2 \sigma^2 )x^2+ \gamma \sigma^2 - 2 \gamma p \sigma x^2 + \frac{1}{2}p^2x^2 \\
        & = -2b \gamma x^4+ \bigl(2\gamma a + \frac{1}{2}(p-2 \gamma \sigma)^2\bigr)x^2+ \gamma \sigma^2.
    \end{align*}
We get {\bf (G)} with any $\gamma > 0$ and this will work for all $p$. We have 
 \begin{align*}
  \Lambda(p) \le \sup_{x \in M} \frac{L_pV(x)}{V(x)} & = \sup_{x \in \R}\Bigl(-2b \gamma x^4+ \bigl(2\gamma a + \frac{1}{2}(p-2 \gamma \sigma)^2\bigr)x^2+ \gamma \sigma^2\Bigr).
  \end{align*}
If $a \ge 0$ choose $\gamma = p/(2\sigma)$.  Then 
 $$
  \Lambda(p) \le   \sup_{x \in \R}\Bigl(-\frac{2b p}{2\sigma}x^4+ \frac{2ap}{2\sigma}x^2+  \frac{p \sigma^2}{2\sigma}\Bigr) 
   = p\left( \frac{a^2}{4b\sigma}+ \frac{\sigma}{2} \right).
   $$
Alternatively, if $a <  0$ choose $\gamma$ so that $2 \sigma \gamma + 2 |a|^{1/2}\gamma^{1/2} = p$. Then 
    $$
    2 \gamma a + \frac{1}{2}(p-2\gamma \sigma)^2 = 2\gamma a+ 2 \gamma |a| = 0
    $$
and $\gamma \sim p/(2\sigma)$, so that 
  $$
  \Lambda(p) \le \sup_{x \in \R}\Bigl(-2b \gamma x^4+ \gamma \sigma^2\Bigr) = \gamma \sigma^2 \sim p\sigma/2
  $$
as $p \to \infty$.  Together
$$
    \limsup_{p \to \infty} \frac{\Lambda(p)}{p} \le  \frac{(\max(a,0))^2}{4b\sigma} + \frac{\sigma}{2}.
    $$    
 
For the lower bound, consider $g(x) = (x^2)^A$ for some $A > 1$.   As in Section \ref{sec corr} we have $L_pg(x) \ge\widetilde{\Gamma}_{p,A}g(x)$ where     
   \begin{align*}
   \widetilde{\Gamma}_{p,A}  
     & = 2A(a-p\sigma) + 2\sqrt{p^2/2-2Ab} \sqrt{\sigma^2 A(2A-1)}
       \end{align*}
so long as $4Ab < p^2$.  For large $p$ choose $A = cp$.  Then for $p > 4bc$ we have  
\begin{align*}
  \Lambda(p) \ge \widetilde{\Gamma}_{p,A} & =2cp(a-p\sigma) + 2\sqrt{p^2-4bcp} \sqrt{\sigma^2 cp(cp-1/2)}\\
   & = 2cap-2c \sigma p^2 + 2c \sigma p^2 \Bigl(1-\frac{2bc}{p} + O(p^{-2}\Bigr)\Bigl(1-\frac{1}{4cp} + O(p^{-2})\Bigr)\\
    & = 2cap + 2c \sigma p \Bigl(-2bc -\frac{1}{4c}\Bigr) + O(1)\\
     & = p\Bigl(-4 \sigma b c^2 +2ac -\frac{\sigma}{2}\Bigr) + O(1)
   \end{align*} 
as $p \to \infty$.  
 If $a >0$ choose $c= a/(4 \sigma b)$.  Then 
    $$
  \liminf _{p \to \infty} \frac{\Lambda(p)}{p} \ge  \frac{a^2}{4b\sigma} - \frac{\sigma}{2}.
     $$ 
Alternatively, if $a \le 0$ choose very small $c > 0$.  Then, in the limit as $c \to 0$ 
    $$
  \liminf _{p \to \infty} \frac{\Lambda(p)}{p} \ge  - \frac{\sigma}{2}.
       $$ 
Together 
$$
    \liminf_{p \to \infty} \frac{\Lambda(p)}{p} \ge \frac{(\max(a,0))^2}{4b\sigma} - \frac{\sigma}{2}.
    $$
In this example the upper and lower estimates do not agree.  However the convexity of $p \mapsto \Lambda(p)$ implies that $\Lambda(p)/p$ is non-decreasing for $p\in \R$, and hence $\lim_{p \to \infty} \Lambda(p)/p$ exists, possibly infinite.  (It exists and is infinite in all the other examples in this section.)  Then the estimates above show 
     $$ 
  \frac{(\max(a,0))^2}{4b\sigma} - \frac{\sigma}{2} \le  \lim_{p \to \infty} \frac{\Lambda(p)}{p}  \le \frac{(\max(a,0))^2}{4b\sigma} + \frac{\sigma}{2}.
   $$
\qed

  \bibliographystyle{plain}
\bibliography{shear.ref}

\end{document}